\documentclass[final,leqno]{siamltex}
\usepackage{amssymb}
\usepackage{amsmath}

\usepackage{graphicx}
\usepackage{multirow}
\usepackage{subfigure}

\newcommand{\bil}[1]{{#1}}

\newtheorem{solver}{Algorithm}[section]
\newtheorem{remark}{Remark}[section]

\newcommand{\hF}{\widehat{F}}
\newcommand{\cP}{\mathcal{P}}
\newcommand{\cQ}{\mathcal{Q}}
\newcommand{\hcP}{\widetilde{\mathcal{P}}}
\newcommand{\cS}{\mathcal{S}}
\newcommand{\cT}{\mathcal{T}}

\newcommand{\bfR}{\mathbf{R}}
\newcommand{\Ome}{\Omega}
\newcommand{\oOme}{\overline{\Omega}}

\begin{document}

\title{Mixed Interior Penalty Discontinuous Galerkin Methods for One-dimensional 
Fully Nonlinear Second Order Elliptic and Parabolic Equations}

\author{Xiaobing Feng\thanks{Department of Mathematics, The University
of Tennessee, Knoxville, TN 37996, U.S.A. (xfeng@math.utk.edu). The work of 
this author was partially supported by the NSF grant DMS-0710831.}
\and
Thomas Lewis\thanks{Department of Mathematics, The University
of Tennessee, Knoxville, TN 37996, U.S.A. (tlewis@math.utk.edu). 
The work of this author was partially supported by the NSF grant DMS-0710831.}
}

\maketitle
\begin{abstract}
This paper is concerned with developing accurate and efficient 
numerical methods for one-dimensional fully nonlinear second order  
elliptic and parabolic partial differential equations (PDEs).
In the paper we present a general framework for constructing high order
interior penalty discontinuous Galerkin (IP-DG) methods for approximating
viscosity solutions of these fully nonlinear PDEs. 
In order to capture discontinuities of the second order derivative $u_{xx}$ 
of the solution $u$, three independent functions $p_1, p_2$
and $p_3$ are introduced to 
represent numerical derivatives using various one-sided limits.
The proposed DG framework, which is based on a nonstandard mixed formulation of
the underlying PDE, embeds a nonlinear problem into a mostly linear system 
of equations where the nonlinearity has been modified to include multiple values 
of the second order derivative $u_{xx}$. 
The proposed framework extends a companion finite difference framework
developed by the authors in \cite{FKL11} and
allows for the approximation of 
fully nonlinear PDEs using high order polynomials and non-uniform meshes. 
In addition to the nonstandard mixed formulation setting, another
main idea is to replace the fully nonlinear differential operator by 
a numerical operator, which is consistent with the differential 
operator and satisfies certain monotonicity (called g-monotonicity) properties. 
To ensure such a g-monotonicity, the crux of the construction 
is to introduce the numerical moment, which plays a critical role in the 
proposed DG framework. The g-monotonicity gives the DG methods
the ability to select the mathematically ``correct" solution
(i.e., the viscosity solution) among all possible solutions.
Moreover, the g-monotonicity allows for the possible development of
more efficient nonlinear solvers as the special nonlinearity
of the algebraic systems can be explored to decouple the equations.
This paper also presents and analyzes numerical results for several 
numerical test problems which are used to guage the accuracy and efficiency 
of the proposed DG methods.
\end{abstract}



\pagestyle{myheadings}
\thispagestyle{plain}
\markboth{XIAOBING FENG AND THOMAS LEWIS}{mIP-DG METHODS FOR
SECOND ORDER FULLY NONLINEAR PDEs}

\section{Introduction}\label{sec-1}

Fully nonlinear partial differential equations (PDEs) refer to a class
nonlinear PDEs which is nonlinear in the highest order derivatives of 
the unknown functions in the equations. Due to their strong nonlinearity,
this class of PDEs are most difficult to analyze analytically and 
to approximate numerically. In the mean time, fully nonlinear PDEs
arise in many applications such as antenna design, astrophysics, 
differential geometry, fluid mechanics, image processing, meteorology, 
mesh generation, optimal control, optimal mass transport, etc \cite{FGN12}, 
 which calls for the development of efficient and reliable numerical 
methods for solving their underlying fully nonlinear PDE problems.

This is the second paper in a series \cite{FKL11} which is 
devoted to developing finite difference (FD) and discontinuous 
Galerkin (DG) methods for approximating {\em viscosity solutions}
of the following general one-dimensional 
fully nonlinear second order elliptic and parabolic equations: 
\begin{equation}\label{pde_ell}
F \left( u_{x x}, u_x, u, x \right) = 0 , \qquad x \in \Omega:=(a,b),
\end{equation}
and
\begin{equation} \label{pde}
u_t + F \left( u_{x x}, u_x, u, t, x \right) = 0 , 
\qquad (x,t) \in \Omega_T:=\Omega\times (0,T),
\end{equation}
which are complemented by appropriate boundary and initial conditions.
The {\em goal} of this paper is to design and implement a class of interior 
penalty discontinuous Galerkin (IP-DG) methods which are based on 
a nonstandard mixed formulation; the proposed IP-DG methods are
named mIP-DG methods. For the ease of presenting the main
ideas and avoiding the technicalities, in this paper we confine 
our attention to the {\em one dimensional} fully nonlinear 
second order PDE problem. The generalization and extension to the high dimensional 
case of the mIP-DG methods of this paper will be presented 
in a forthcoming work. In fact, it will be seen 
later that even in the one dimensional case, the construction 
and analysis of the proposed mIP-DG methods are already quite complicated.  

It is well known \cite{FGN12} that the primary challenges 
for approximating viscosity solutions of fully nonlinear PDEs 
are caused by the very notion of viscosity solutions themselves (see section \ref{sec-2}
for the definition). Unlike the notion of weak solutions for linear and 
quasilinear PDEs, the notion of viscosity solutions by design is 
non-variational, and, in general, viscosity solutions do not satisfy the 
underlying PDEs in a tangible sense. The non-variational nature of viscosity 
solutions immediately prevents any attempt to directly and 
straightforwardly construct Galerkin-type (including DG) methods 
for approximating fully nonlinear PDEs; in other words, nonlinearity in the highest 
order derivatives of the unknown function does not allow 
one to perform integration by parts to transfer one order of derivatives
to test functions as often done with linear and quasilinear PDEs.     
Another big challenge for approximating viscosity solutions of 
fully nonlinear PDEs is caused by the conditional uniqueness of
viscosity solutions, namely, viscosity solutions 
may only be unique in a restricted function class. Requiring 
numerical solutions to stay or approximately stay in the same function 
class often imposes a difficult constraint for designing numerical methods.
Finally, we like to mention that as expected, solving the resulting 
strong nonlinear (algebraic) systems, regardless which discretization 
method is used, is another difficult issue encountered with numerical 
fully nonlinear PDEs.

The mIP-DG methods proposed in this paper aim to approximate 
viscosity solutions of \eqref{pde_ell} and \eqref{pde} which 
belong to $H^1(\Omega)$ in the spatial variable. We note that such 
a viscosity still does not satisfy the underlying PDEs in 
a tangible sense. We also mention that in order to approximate viscosity 
solutions that do not have $H^1$ regularity in the spatial variable, 
we refer the reader to
a companion paper \cite{Feng_Lewis12c} in which we propose another class of more 
complicated mixed discontinuous Galerkin that incorporates a local 
discontinuous Galerkin (LDG) approach instead of the IP-DG approach. 
Such an alternate LDG approach is also more appropriate when a more 
accurate approximation for $u_x$ is desired. 

Several novel ideas are utilized to design the mIP-DG methods 
in this paper which are briefly described below. 
Since integration by parts, which is the necessary tool for
constructing any DG method, cannot be performed on equation \eqref{pde_ell}, 
{\em the first key idea} is to introduce the auxiliary variable $p:=u_{xx}$ and 
rewrite the original fully nonlinear PDE as a system of PDEs:  
\begin{eqnarray} \label{mixed_1}
F(p,u_x,u,x) &=0, \\
p-u_{xx} &=0. \label{mixed_2}
\end{eqnarray}
Unfortunately, since $u_{xx}$ may not exist for a viscosity 
solution $u\in H^1(\Omega)$, the the above mixed formulation may 
not make sense. To overcome this difficulty, {\em the second key idea} 
is to replace $p:=u_{xx}$ by three possible values of 
$u_{xx}$, namely, the left and right limits, as well as their average.  
Thus, we have
\begin{eqnarray} \label{p1}
p_1(x)-u_{xx}(x^-) &=0, \\
p_2(x)-u_{xx}(x^a) &=0, \label{p2} \\
p_3(x)-u_{xx}(x^+) &=0, \label{p3}
\end{eqnarray}
where $u_{xx}(x^a)$ can be thought of as the arithmetic average of $u_{xx}(x^-)$
and $u_{xx}(x^+)$. We remark that \eqref{p1} and \eqref{p3} can be regarded 
as two ``one-sided" Poisson problems in $u$, and \eqref{p2} can be thought of as
the ``regular" Poisson problem. To incorporate the multiple values of $u_{xx}$, 
equation \eqref{mixed_1} must be modified because $F$ is only defined for a 
single value function $p$. To this end, we need 
{\em the third key idea} of this paper, which is to replace \eqref{mixed_1} by 
\begin{eqnarray} \label{F_hat}
\hF(p_1, p_2, p_3, u_x,u,x) =0,
\end{eqnarray} 
where $\hF$, which is called a {\em numerical operator}, should be some
well-chosen approximation to $F$. 

Natural questions that now arise are what are the criterions for $\hF$, and
how can such a numerical operator $\hF$ be constucted? 
These are two immediate questions which must be 
addressed. To do so, we need {\em the fourth key idea} of this paper, 
which is to borrow and adapt the notion of the numerical 
operators from our previous work \cite{FKL11} where a general
finite difference framework has been developed for fully nonlinear 
second order PDEs. In summary, the criterions for $\hF$ include 
{\em consistency} and {\em g-monotonicity} (generalized monotonicity), 
for which precise definitions can be found in section \ref{sec-2}.  
It should be 
pointed out that in order to construct the desired numerical operator
$\hF$, a fundamental idea used in \cite{FKL11} is to introduce the 
concept of {\em the numerical moment}, which can be regarded as 
a direct numerical realization for the moment term in {\em the vanishing
moment methodology} introduced in \cite{Feng_Neilan08} (also see
\cite[section 4]{FGN12},\cite{Feng_Neilan11}). Finally, we need to design 
a DG discretization for the mixed system \eqref{p1}--\eqref{F_hat} 
to accomplish the goal. {\em The fifth key idea} of this paper is
to use different {\em numerical fluxes} in the formulations of
IP-DG methods for the ``one-side" Poisson problems \eqref{p1} and \eqref{p3}
as well as for the ``regular" Poisson problem \eqref{p2}. We remark that,
to the best of our knowledge, this is one of a few scenarios 
in numerical PDEs where the flexibility and superiority 
(over other numerical methodologies)
of the DG methodology makes a vital difference.
 
The remainder of this paper is organized as follows. In section \ref{sec-2}
we collect some preliminaries including the definition of viscosity 
solutions, the definitions of the consistency and g-monotonicity of numerical 
operators, and the concept of the numerical moment.  
In section \ref{sec-3} we present the detailed formulation 
of mIP-DG methods for fully 
nonlinear elliptic equation \eqref{pde_ell} following the outline 
described above. In section \ref{sec-4} we consider both explicit and 
implicit in time fully discrete mIP-DG methods for fully nonlinear 
parabolic equation \eqref{pde} based on the method of lines approach. 
The forward and backward Euler time-stepping schemes combined with 
the spatial mIP-DG methods will be specifically formulated. 
In section \ref{sec-5} we present many numerical experiments 
for the proposed mIP-DG methods and their fully discrete 
counterparts for the parabolic equation \eqref{pde}. 
These numerical experiments verify the accuracy of the proposed 
mIP-DG methods and also demonstrate the efficiency of these methods.
Finally, we complete the paper with a brief summary and 
some concluding remarks in section \ref{sec-6}.

\section{Preliminaries} \label{sec-2}
For a bounded open domain $\Ome\subset\mathbf{R}^d$, let $B(\Ome)$, 
$USC(\Ome)$ and $LSC(\Ome)$ denote, respectively, the spaces of bounded,
upper semi-continuous, and lower semicontinuous functions on $\Ome$.
For any $v\in B(\Ome)$, we define
\[
v^*(x):=\limsup_{y\to x} v(y) \qquad\mbox{and}\qquad
v_*(x):=\liminf_{y\to x} v(y). 
\]
Then, $v^*\in USC(\Ome)$ and $v_*\in LSC(\Ome)$, and they are called
{\em the upper and lower semicontinuous envelopes} of $v$, respectively.

Given a bounded function $F: \cS^{d\times d}\times\mathbf{R}^d\times 
\mathbf{R}\times \oOme \to \mathbf{R}$, where $\cS^{d\times d}$ denotes the set
of $d\times d$ symmetric real matrices, the general second order
fully nonlinear PDE takes the form
\begin{align}\label{e2.1}
F(D^2u,\nabla u, u, x) = 0 \qquad\mbox{in } \oOme.
\end{align}
Note that here we have used the convention of writing the boundary condition as a
discontinuity of the PDE (cf. \cite[p.274]{Barles_Souganidis91}).

The following two definitions can be found in \cite{Gilbarg_Trudinger01,
Caffarelli_Cabre95,Barles_Souganidis91}.

\begin{definition}\label{def2.1}
Equation \eqref{e2.1} is said to be elliptic if for all
$(\mathbf{q},\lambda,x)\in \mathbf{R}^d\times \mathbf{R}\times \oOme$ there holds
\begin{align}\label{e2.2}
F(A, \mathbf{q}, \lambda, x) \leq F(B, \mathbf{q}, \lambda, x) \qquad\forall 
A,B\in \cS^{d\times d},\, A\geq B, 
\end{align}
where $A\geq B$ means that $A-B$ is a nonnegative definite matrix.
\end{definition}
We note that when $F$ is differentiable, the ellipticity
also can be defined by requiring that the matrix $\frac{\partial F}{\partial A}$
is negative semi-definite (cf. \cite[p. 441]{Gilbarg_Trudinger01}).

\begin{definition}\label{def2.2}
A function $u\in B(\Ome)$ is called a viscosity subsolution (resp.
supersolution) of \eqref{e2.1} if, for all $\varphi\in C^2(\oOme)$,
if $u^*-\varphi$ (resp. $u_*-\varphi$) has a local maximum
(resp. minimum) at $x_0\in \oOme$, then we have
\[
F_*(D^2\varphi(x_0),\nabla \varphi(x_0), u^*(x_0), x_0) \leq 0 
\]
(resp. $F^*(D^2\varphi(x_0),\nabla \varphi(x_0), u_*(x_0), x_0) \geq 0$).
The function $u$ is said to be a viscosity solution of \eqref{e2.1}
if it is simultaneously a viscosity subsolution and a viscosity
supersolution of \eqref{e2.1}.
\end{definition}


We remark that if $F$ and $u$ are continuous, then the upper and lower $*$
indices can be removed in Definition \ref{def2.2}. The definition
of ellipticity implies that the differential operator $F$
must be non-increasing in its first argument in order to be
elliptic. It turns out that ellipticity provides a sufficient
condition for equation \eqref{e2.1} to fulfill a maximum principle
(cf. \cite{Gilbarg_Trudinger01, Caffarelli_Cabre95}).
It is clear from the above definition that viscosity solutions
in general do not satisfy the underlying PDEs in a tangible sense, and
the concept of viscosity solutions is {\em nonvariational}. Such
a solution is not defined through integration by parts against arbitrary test
functions; hence, it does not satisfy an integral identity. As pointed
out in section \ref{sec-1}, the nonvariational nature of viscosity
solutions is the main obstacle that prevents direct construction
of Galerkin-type methods, which are based on variational formulations.

\smallskip
The following definitions are adapted from \cite{FKL11} in the case $d=1$.

\smallskip
\begin{definition}\label{def2.3}
\begin{itemize}
\item[{\rm (i)}] A function $\hF: \bfR^6\to \bfR$ is called  
a {\em numerical operator}. 
\item[{\rm (ii)}] A numerical operator $\hF$ is said to be consistent (with 
the differential operator $F$) if $\hF$ satisfies
\begin{align}\label{A1a}
\liminf_{p_k\to p, k=1,2,3\atop q_1\to q, \lambda_1\to \lambda,\xi_1\to \xi} 
\hF(p_1,p_2,p_3,q_1,\lambda_1, \xi_1) \geq F_*(p,q,\lambda,\xi),\\
\limsup_{p_k\to p, k=1,2,3 \atop q_1\to q, \lambda_1\to \lambda,\xi_1\to \xi} 
\hF(p_1,p_2,p_3, q_1, \lambda_1,\xi_1) 
\leq F^*(p,q,\lambda,\xi), \label{A1b}
\end{align}
where $F_*$ and $F^*$ denote respectively the lower and the upper
semi-continuous envelopes of $F$.

\item[{\rm (iii)}] A numerical operator $\hF$ is said to be
{\em g-monotone} if $\hF(p_1,p_2,p_3,q,\lambda,\xi)$ is monotone increasing 
in $p_1$ and $p_3$ and monotone decreasing in $p_2$, that is, 
$\hF(\uparrow,\downarrow,\uparrow,q,\lambda,\xi)$.
\end{itemize}

\end{definition}

\medskip
We note that the above consistency and g-monotonicity play a critical role in 
the finite difference framework established in \cite{FKL11}. They also play
an equally critical role in the mIP-DG methods of this paper. 
We also note that in practice the consistency is easy to fulfill and 
to verify, but the g-monotonicity is not. In order to ensure the 
g-monotonicity, one key idea of \cite{FKL11} is to introduce the 
concept of {\em the numerical moment} to help.  The following 
are two examples of so-called Lax-Friedrichs-like numerical 
operators \cite{FKL11}:
\begin{align}\label{LF1}
\hF_1(p_1,p_2,p_3,q,\lambda,\xi)&:= F(p_2,q,\lambda,\xi) 
+ \alpha_1\bigl(p_1-2p_2+p_3\bigr),\\
\hF_2(p_1,p_2,p_3,q,\lambda,\xi)&:=F\Bigl(\frac{p_1+p_2+p_3}{3},q,\lambda,\xi \Bigr) 
+ \alpha_2\bigl(p_1-2p_2+p_3\bigr),\label{LF2}
\end{align}
where $\alpha_1$ and $\alpha_2$ are undetermined positive constants and 
the last term in \eqref{LF1} and \eqref{LF2} is called {\em the numerical moment}. 
It is trivial to verify that $\hF_1$ and $\hF_2$ are consistent with $F$.
To ensure $\hF_1$ to be g-monotone,  we need 
\begin{equation} \label{alpha_ell}
\alpha > \left| \frac{\partial F}{\partial u^{\prime \prime}} \right|, 
\end{equation}
assuming adequate regularity for the operator $F$. We remark that
it is natural to require that $\hF_1$ is decreasing in $p_2$ because
by the definition of ellipticity, $F$ is decreasing in $u^{\prime \prime}$.

\section{Formulation of mIP-DG methods for elliptic problems}\label{sec-3}
 
We first consider the elliptic problem $(\ref{pde_ell})$ with boundary conditions
\begin{equation}
u(a) = u_a \quad\mbox{and}\quad  u(b) = u_b \label{bc_ell} 
\end{equation}
for two given constants $u_a$ and $u_b$. 

Let $\left\{ x_j \right\}_{j=0}^J \subset \overline{\Omega}$ be a mesh 
for $\overline{\Omega}$ such that $x_0 = a$ and $x_J = b$.  
Define $I_{j} = \left( x_{j-1} , x_{j} \right)$ and $h_j = x_{j} - x_{j-1}$
for all $j = 1, 2, \ldots, J$, $h_0 = h_{J+1} = 0$ and 
$h = \max_{1 \leq j \leq J} h_j$.  Let $\mathcal{T}_h$ denote
the collection of the intervals $\{I_j\}_{j=1}^J$ which form a partition 
of the domain $\overline{\Omega}$. We also introduce the broken $H^1$-space 
\[
H^1(\cT_h):= \prod_{I\in \cT_h} H^1(I)
\]
and the broken $L^2$-inner product
\[
(v ,w)_{\mathcal{T}_h}:= \sum_{j=1}^J \int_{I_j} v w\, dx \qquad 
\forall v,w\in H^1(\cT_h).
\]
For a fixed integer $r \geq 1$, we define the standard DG finite element space
$V^h \subset H^1(\cT_h)\subset L^2 (\mathcal{T}_h)$ by
\[
V^h := \prod_{I \in \mathcal{T}_h} \mathcal{P}_{r} (I),
\]
where $ \mathcal{P}_{r} (I)$ denotes the set of all polynomials on $I$ with
degree not exceeding $r$.  We also introduce the following
standard jump and average notations:
\begin{align*}
[v_h(x_j)] &:= v_h(x_j^-)-v_h(x_j^+) \qquad\mbox{for } j=1,2,\cdots, J-1,\\
[v_h(x_0)] &:=-v_h(x_0), \qquad [v_h(x_J)] := v_h(x_J);\\
\{v_h(x_j)\} &:=\frac12\Bigl( v_h(x_j^-)+v_h(x_j^+) \Bigr) \qquad\mbox{for } 
j=1,2,\cdots, J-1,\\
\{v_h(x_0)\} &:=v_h(x_0), \qquad \{v_h(x_J)\} := v_h(x_J).
\end{align*}
It is trivial to verify the following so-called ``magic formulas":
\begin{align} \label{jump1}
[v(x_j) w(x_j)] 
&= v ( x_j^-) [w(x_j)] + [v(x_j)] w(x_j^+),\\
[v(x_j) w(x_j)] 
&=\{v (x_j)\}[w(x_j)] + [v(x_j)]\{ w(x_j)\},  \label{jump2} \\
[v(x_j) w(x_j)] 
&=v(x_j^+) [w(x_j)] + [v(x_j)] w(x_j^-).  \label{jump3}  
\end{align}

Let $\gamma_{0i} > 0$ for $i = 1, 2, 3$ denote interior penalty 
parameters. It will be clear later that to avoid redundancy 
of three equations for $p_1, p_2$ and $p_3$, we need to require that
$\gamma_{02} > \max \left\{ \gamma_{01} , \gamma_{02} \right\}$.
Define the interior penalty terms
\begin{equation} \label{j0}
J_{0i} \left( v, w \right) = \sum_{j = 0}^{J} \frac{\gamma_{0i}}{h_{j, j+1}} 
\left[ v \left( x_{j} \right) \right] \; \left[ w \left( x_{j} \right) \right]  , 
\end{equation}
for $i = 1, 2, 3$, where
\[
h_{j, j+1} = \max \left\{ h_j, h_{j+1} \right\} \qquad\mbox{for }
j = 0,1, 2, \ldots, J.
\]

We now are ready to formulate our DG discretizations for equations 
\eqref{p1}--\eqref{F_hat}. First, for (fully) nonlinear equation \eqref{F_hat}
we simply approximate it by its broken $L^2$-projection into $V^h$, namely,
\begin{equation}\label{pde_ell_weak}
\bil{a}_0 \bigl(u_h , p_{1h}, p_{2h}, p_{3h}; \phi_{0h} \bigr)  = 0 
\qquad \forall \phi_{0h} \in V^h, 
\end{equation}
where
\[
\bil{a}_0 (u, p_1, p_2, p_3; \phi_{0}) 
=\Bigl(\hF(p_1, p_2, p_3, u^\prime, u,\cdot),\phi_{0} \Bigr)_{\mathcal{T}_h}. 
\]

Next, we discretize the three {\em linear} equations \eqref{p1}--\eqref{p3}.
Notice that for given ``sources" $\{p_i\}_{i=1}^3$, \eqref{p1}--\eqref{p3}
are three (different) Poisson equations for $u$. Thus, we can use the 
standard IP-DG formulation for the Laplacian operator to discretize these
equations. However, there is a crucial distinction for doing so on the three
equations, that is, we use, respectively, ``magic formulas" \eqref{jump1}, 
\eqref{jump2}, and \eqref{jump3} when we add the local integration by parts
formula to handle the jump terms at the interior nodes. To realize the
above strategy, we define the bilinear forms 
$\bil{b}_i: H^1(\mathcal{T}_h)\times H^1(\mathcal{T}_h)\to \mathbb{R}$ by 
\begin{align} \label{blin}
\bil{b}_i (v, w)
&:=(v^{\prime}, w^{\prime})_{\mathcal{T}_h} 
+ v^{\prime}(a) w(a) - \epsilon \, v(a) w^{\prime}(a) - v^{\prime} (b) w(b) \\
\nonumber & \qquad	
+ \epsilon \, v(b) w^{\prime} (b) 
+ J_{0i} (v, w) \quad\forall v, w \in H^1 \left( \mathcal{T}_h \right), \,
i = 1, 2, 3,
\end{align}
where $\epsilon \in \{-1,0,1\}$ is often called the ``symmetrization"
parameter \cite{Riviere08}.  Note,  $\bil{b}_i$ is symmetric if $\epsilon=-1$, 
nonsymmetric if $\epsilon=1$, and incomplete if $\epsilon=0$. 
Using the bilinear forms $\bil{b}_i$, we define the following 
DG discretizations of \eqref{p1}--\eqref{p3}:
\begin{align} \label{pi_bil}
\bil{a}_i (u_h, p_{ih};\phi_{ih} ) 
& = \bil{f}_i (\phi_{ih}), \qquad 
\forall \phi_{ih} \in V^h, \,\, i=1,2,3, 
\end{align}
where 
\begin{align*}
&\bil{a}_1 (u, p_1 ; \phi_{1})
=( p_1 , \phi_{1})_{\Omega} + \bil{b}_i (u, \phi_{1}) 
-\sum_{j=1}^{J-1} \Bigl( u^{\prime} (x_{j}^-) \bigl[ \phi_{1} (x_{j}) \bigr]
- \epsilon \, \bigl[u( x_{j}) \bigr] \phi_{1}^{\prime} (x_{j}^-) \Bigr), \\
&\bil{a}_2 (u, p_2; \phi_{2} )
= (p_2, \phi_{2} )_{\Omega} + \bil{b}_i (u, \phi_{2} ) 
-\sum_{j=1}^{J-1}\Bigl( \bigl\{u^{\prime}(x_{j})\bigr\}\bigl[\phi_{2}(x_{j})\bigr]
- \epsilon \, \bigl[u(x_{j})\bigr] \bigl\{\phi_{2}^{\prime}(x_{j})\bigr\} \Bigr),\\
&\bil{a}_3 (u, p_3; \phi_{3})  
= (p_3, \phi_{3})_{\Omega} + \bil{b}_i (u,\phi_{3})
-\sum_{j=1}^{J-1} \Bigl( u^{\prime} (x_{j}^+) \bigl[\phi_{3}( x_{j}) \bigr]
-\epsilon \, \bigl[u( x_{j})\bigr] \phi_{3}^{\prime}(x_{j}^+) \Bigr), 
\end{align*}
and
\[
\bil{f}_i (\phi_{i}) = \Bigl(\frac{\gamma_{0i}}{h_{0, 1}} \; \phi_{i}(a) 
-\epsilon\, \phi^{\prime}_{i}(a) \Bigr) u_a 
+\Bigl( \frac{\gamma_{0i}}{h_{J, J+1}} \; \phi_{i} (b) 
+\epsilon \, \phi^{\prime}_{i}(b) \Bigr) u_b  
\]
for $i = 1, 2, 3$. 

In summary, our mIP-DG methods for the fully nonlinear Dirichlet problem 
\eqref{pde_ell}, \eqref{pde}, and \eqref{bc_ell} are defined as seeking
$\bigl(u_h,p_{1h},p_{2h},p_{3h} \bigr)\in [V^h]^4$ such that \eqref{pde_ell_weak} 
and \eqref{pi_bil} hold.

We conclude this section with a few remarks.

\smallskip
\begin{remark}
(a) Looking backwards, \eqref{pi_bil} provides a proper interpretation for 
each of $p_{1h}$, $p_{2h}$, and $p_{3h}$ for a given function $u_h$. Each 
$p_{ih}$ defines a discrete second order derivative of $u_h$. 
The functions $p_{1h}$, $p_{2h}$, and $p_{3h}$ should be very close to each other 
if $u_{xx}$ exists; however, their discrepancies are expected to be large 
if $u_{xx}$ does not exist. $p_{1h}$, $p_{2h}$, and $p_{3h}$ defined by \eqref{pi_bil}
can be regarded as high order extensions of their lower order counterparts 
defined in \cite{FKL11}.

(b) It is easy to check that the three equations defined by \eqref{pi_bil}
are linearly independent provided that $\gamma_{02} > \max \left\{ \gamma_{01}, 
\gamma_{02} \right\}$.

(c) The reason for $r \neq 0$ can be explained as follows. When $r=0$, 
the piecewise constant basis functions have piecewise zero derivatives 
on the given mesh. After eliminating the jump terms containing derivatives
in \eqref{pi_bil}, it is clear that the ability for $p_1$ and $p_3$ 
to carry information from the left and the right, respectively, is lost.
Furthermore, if $\gamma_{01} = \gamma_{02} = \gamma_{03}$, 
then $\bil{a}_1 = \bil{a}_2 = \bil{a}_3$, which in turn implies that
on a uniform mesh $p_{1h} = p_{2h} = p_{3h}$, they all are equal to the centered
difference approximation for the second order derivative of $u_h$.
As a result, the numerical moment term vanishes and we are left
with the standard three-point finite difference approximation for \eqref{pde_ell}
and \eqref{bc_ell}. 
On the other hand, when $r \geq 1$, the numerical operator
maintains the directional interpretations for $p_1$ and $p_3$, allowing
the numerical operator to take advantage of the numerical moment.  


(e) Notice that \eqref{pde_ell_weak}--\eqref{pi_bil} is a nonlinear 
system of equations, with the nonlinearity only appearing in $\bil{a}_0$. 
Thus, a nonlinear solver is necessary in implementing the above scheme. 
In section \ref{sec-5}, 
an iterative method is used with initial guess given by projecting 
the secant line resulting from the boundary conditions into $V^h$.  
Since a good initial guess is essential for most nonlinear solvers 
to converge, another possibility is to first linearize the nonlinear 
operator and solve the resulting linear system first.  
However, we show in our numerical tests that the simple initial 
guess works well in many cases. We suspect that the g-monotonicity 
of $\hF$ enlarges the domain of ``good" initial values over which 
the iterative method converges.
\end{remark}

\section{Formulation of fully discrete mIP-DG methods for parabolic problems}
\label{sec-4}

The goal of this section is to extend our mIP-DG methods to solving 
the initial-boundary value problem for \eqref{pde} using the method of lines.  
Let the initial condition be given by
\begin{equation} \label{ic}
u(0, x) = u_0 (x) \qquad\forall x\in \Omega
\end{equation}
and the boundary conditions be given by
\begin{equation} \label{bc}
u(t,a) = u_a(t), \qquad  u(t,b) = u_b (t) \qquad \forall  t\in (0,T)..
\end{equation}
For the ease of the presentation, we shall only consider the forward 
and backward Euler methods, although higher order time-stepping methods 
can also be formulated. 

For a fixed integer $M > 0$ and let $\Delta t = \frac{T}{M}$ be the 
time-step size.  Define $t^n:=n\Delta t$ for $n = 0, 1, \cdots, M$. 
Then, the forward Euler method (in operator form) for \eqref{pde} is
defined by seeking $u^{n+1}: \Ome \to \bfR$ such that
\begin{equation}\label{e4.1} 
u^{n+1} = u^{n} - 
\Delta t \, F \left( u_{x x}^n , u_x^n , u^n , t^n, \cdot \right) 
\qquad \mbox{in } \Omega
\end{equation}
and the standard backward Euler method (in operator form) for \eqref{pde} is 
defined by seeking $u^{n+1}: \Ome\to \mathbb{R}$ such that
\begin{equation}\label{e4.2}
u^{n+1} = u^{n} - \Delta t \, F \bigl(u_{x x}^{n+1}, u_x^{n+1}, u^{n+1}, t^{n+1}, 
\cdot \bigr) \qquad\mbox{in } \Omega
\end{equation}
for $n = 0, 1, \ldots, M-1$, where
\[
u^0 = u_0 \qquad \mbox{in } \Omega.
\]
The compatibility condition between the initial and boundary data 
immediately implies that
\begin{equation}\label{e4.3}
u^n (a) = u_a(t^n), \qquad u^n (b) = u_b(t^n), \quad n = 1, 2, \cdots, M.
\end{equation}

We next apply the mIP-DG framework developed in the previous section to 
equations \eqref{e4.1} and \eqref{e4.2} for their spatial discretizations.
We present these two cases separately below because they require 
different treatments and involve different technicalities.
 
\subsection{Forward Euler method} \label{sec-4.1}
It turns out that the forward Euler method is tricky to formulate because 
the variables $u_h^{n+1}$ and $p_{jh}^{n+1}$ ($j=1,2,3$) are not determined 
simultaneously.  Instead, they are constructed sequentially.  
For a given $u_h^n$, we first construct 
$p_{jh}^n$ for $j=1,2,3$ using \eqref{pi_bil}.  We then 
define $u_h^{n+1}$ to be a modified $L^2$-projection of the right-hand side 
of \eqref{e4.1}. To take care of the boundary condition, we choose to 
enforce the boundary condition for $u_h^{n+1}$ weakly in the 
definition of the modified $L^2$-projection. 

Specifically, for any $v\in L^2(\Ome)$, we recall that the standard 
$L^2$-projection $\cP_h v\in V^h$ of $v$ is defined by 
\begin{equation} \label{L2_proj}
\bigl( \cP_h v, \phi_h \bigr)_{\cT_h} = \bigl(v, \phi_h \bigr)_{\cT_h}
\qquad\forall \phi_h\in V^h.
\end{equation}
For any $v\in C^0(\overline{\Ome})$, we define a
modified $L^2$-projection $\hcP_h v\in V^h$ of $v$ by
\begin{align} \label{mL2_proj}
\bigl( \hcP_h v, \phi_h \bigr)_{\cT_h} 
&+ \frac{1}{\sqrt{h}} \Bigl(\hcP_h v(a) \phi_h(a) + \hcP_h v(b) \phi_h(b) \Bigr) \\
&\quad = \bigl(v, \phi_h \bigr)_{\cT_h} 
+ \frac{1}{\sqrt{h}} \Bigl(v(a) \phi_h(a) + v(b) \phi_h(b) \Bigr)
\qquad\forall \phi_h\in V^h, \nonumber
\end{align}
and the corresponding modified $L^2$-projection operator $\hcP_h: 
L^2(\Ome)\cap C^0(\overline{\Ome})\to V^h$. In the above definition, the 
boundary condition \eqref{bc} is weakly enforced via a penalty technique which 
is due to Nitsche \cite{Nitsche70}.  
	
For a given function $v\in H^1(\cT_h)$ which satisfies the boundary 
condition \eqref{bc}, we define its three discrete ``one-side" second order 
derivatives $\cQ_h^- v, \cQ_h^a v, \cQ_h^+ v\in V^h$ ($j=1,2,3$) 
using \eqref{pi_bil} as follows:
\begin{align} \nonumber 
\bigl( \cQ_h^- v,\phi_h \bigr)_\Omega
& =\Bigl(\frac{\gamma_{01}}{h_{0,1}}\phi_h(a) 
-\epsilon\phi_h^{\prime}(a)\Bigr) u_a(t)
+\Bigl(\frac{\gamma_{01}}{h_{J, J+1}}\phi_h(b)
+\epsilon\phi_h^{\prime}(b) \Bigr) u_b(t) \\
&\quad -\bil{b}_1(v,\phi_h)
+\sum_{j=1}^{J-1} \Bigl( v^{\prime}(x_{j}^-) [\phi_h(x_{j})]
-\epsilon [v(x_{j})] \phi_h^{\prime}(x_{j}^-) \Bigr) 
\quad\forall \phi_h\in V^h, \label{projection_p1} \\
\bigl(\cQ_h^a v,\phi_h \bigr)_\Omega
&=\Bigl( \frac{\gamma_{02}}{h_{0,1}} \phi_2(a) 
-\epsilon \phi_h^{\prime}(a) \Bigr) u_a(t)
+\Bigl( \frac{\gamma_{02}}{h_{J, J+1}} \phi_h(b)
+\epsilon \phi_h^{\prime}(b)  \Bigr) u_b (t) \nonumber \\
& \quad -\bil{b}_2(v, \phi_h) 
+\sum_{j=1}^{J-1} \Bigl(\{v^{\prime}(x_{j})\} [\phi_h(x_{j})]
-\epsilon [v(x_{j})] \{\phi_h^{\prime}(x_{j})\} \Bigr)  
\quad\forall \phi_h\in V^h, \label{projection_p2} \\
\Bigl( \cQ_h^+ v,\phi_h \bigr)_\Omega
& = \Bigl( \frac{\gamma_{03}}{h_{0,1}} \phi_h(a) 
-\epsilon \phi_h^{\prime}(a) \Bigr) u_a(t)
+ \Bigl( \frac{\gamma_{03}}{h_{J,J+1}} \phi_h(b)
+ \epsilon \phi_h^{\prime}(b)  \Bigr) u_b (t) \nonumber \\
& \quad -\bil{b}_3( v,\phi_h) 
+ \sum_{j=1}^{J-1} \Bigl( v^{\prime}(x_{j}^+) [\phi_3(x_{j})]
- \epsilon [v(x_{j})] \phi_h^{\prime}(x_{j}^+) \Bigr) 
\quad\forall \phi_h\in V^h, \label{projection_p3} 
\end{align}
and the corresponding operators $\cQ_h^-, \cQ_h^a, \cQ_h^+: H^1(\cT_h)\to V^h$.

With the help of the operators $\cQ_h^-, \cQ_h^a, \cQ_h^+$, and $\hcP_h$, 
we now define our fully discrete forward Euler method for the initial-boundary
value problem \eqref{pde}, \eqref{ic}, \eqref{bc} as follows:
for $n = 0, 1,\ldots, M-1$,
\begin{align} \label{forward_euler}
u_h^{n+1} &=\hcP_h \Bigl( u_h^{n} - \Delta t \, \widehat{F} 
\bigl( \cQ_h^- u_h^n, \cQ_h^a u_h^n ,\cQ_h^+ u_h^n,
\{u_{hx}^n\}, \{u_h^n\}, t^n, \cdot \bigr) \Bigr), \\
u_h^0 &= \cP_h u_0, \label{forward_euler_1}
\end{align}
where $\{u_{hx}^n(x)\}$ and $\{u_h^n(x)\}$ denote the average values of 
$u_{hx}^n$ and $u_h^n$ at $x$.  We remark that the average operation is 
needed because the approximation $u_h^n$ is only defined piecewisely on the mesh
$\cT_h$. Moreover, a projection is necessary at each time step 
since $\hF$ may not belong to $V^h$ and the prescribed Dirichlet boundary 
condition must be enforced at each time step.  

\subsection{Backward Euler method} \label{sec-4.2}
We rewrite the backward Euler scheme \eqref{e4.2} as
\begin{equation}\label{e4.13}
u^{n} + \Delta t \, F\bigl(u_{xx}^{n}, u_x^{n}, u^{n}, t^{n}, 
\cdot\bigr) = u^{n-1}
\end{equation}
for $n=1,\cdots,M$, where $u^0 =u_0$. 

Clearly, \eqref{e4.13} and \eqref{pde} have the same form. Thus, the
spatial discretization of the backward Euler scheme \eqref{e4.2} is 
a straightforward adaptation of the mIP-DG framework for elliptic PDEs
developed in section \ref{sec-3}. To the end, we define a new 
numerical operator $\widehat{G}$ by 
\begin{equation} \label{Ghat}
\widehat{G} \bigl(p_1, p_2, p_3, u, t, x \bigr) 
:=u(x,t)+ \Delta t \, \hF\bigl(p_1, p_2, p_3, u_x, u, t, x\bigr)
\qquad\forall (x,t)\in \Ome_T.
\end{equation}
Then, analogous to the formulation for \eqref{pde_ell} and \eqref{bc_ell}, 
the fully discrete backward Euler mIP-DG methods for \eqref{pde}, \eqref{ic}, 
and \eqref{bc} is defined by seeking $(u_h^n, p_{1h}^n, p_{2h}^n, p_{3h}^n)
\in (V^h)^4$ such that, for $n = 1, 2, \ldots, M$, 
\begin{align} \label{backward_euler}
\bil{\widehat{a}}_0 \bigl(t^n, u_h^n, p_{1h}^n, p_{2h}^n, p_{3h}^n; \phi_{0h} \bigr)
&= \bigl(u_h^{n-1}, \phi_{0h} \bigr)_{\mathcal{T}_h}
\qquad \forall\phi_{0h} \in V^h \\
\bil{\widehat{a}}_i \bigl(u_h^n, p_{ih}^n; \phi_{ih} \bigr)
&= \bil{g}_i(t^n, \phi_{ih}) \qquad\forall \phi_{ih} \in V^h, \,\, i=1,2,3,\\
u_h^0 &= \cP_h u_0. \label{backward_euler_1} 
\end{align}
where
\begin{align*}
&\bil{\widehat{a}}_0 \bigl(t, u, p_1, p_2, p_3; \phi_0 \bigr)
=\Bigl(\widehat{G} \bigl( p_1, p_2, p_3 , u , t, \cdot \bigr), \phi_{0} 
\Bigr)_{\mathcal{T}_h} \\
&\bil{\widehat{a}}_1 \bigl(u, p_1; \phi_{1}\bigr)
=( p_1, \phi_{1})_{\Omega} + \bil{b}_1 (u, \phi_{1}) 
- \sum_{j=1}^{J-1} \Bigl( u^{\prime}( x_{j}^-) [\phi_{1}( x_{j})]
- \epsilon \, [u(x_{j})] \phi_{1}^{\prime}( x_{j}^-) \Bigr), \\
&\bil{\widehat{a}}_2 \bigl(u, p_2; \phi_{2} \bigr)
=( p_2, \phi_{2})_{\Omega} + \bil{b}_2 (u, \phi_{2} ) 
-\sum_{j=1}^{J-1} \Bigl(\{u^{\prime} (x_{j})\} [\phi_{2}(x_{j})] 
-\epsilon \, [u(x_{j})] \{\phi_{2}^{\prime}(x_{j})\} \Bigr), \\
&\bil{\widehat{a}}_3 \bigl(u, p_3; \phi_{3} \bigr)
 = ( p_3, \phi_{3} )_{\Omega} + \bil{b}_3 (u, \phi_{3}) 
- \sum_{j=1}^{J-1} \Bigl( u^{\prime}(x_{j}^+) [\phi_{3}( x_j)]
- \epsilon \, [u(x_{j})] \phi_{3}^{\prime}( x_{j}^+ ) \Bigr), 
\end{align*}
and
\begin{align*}
\bil{g}_i (t, \phi_{i} )
& = \Bigl( \frac{\gamma_{0i}}{h_{0, 1}} \; \phi_{i}(a) 
- \epsilon \, \phi^{\prime}_{i}(a) \Bigr) u_a(t)
+ \Bigl( \frac{\gamma_{0i}}{h_{J, J+1}} \; \phi_{i}(b)
+ \epsilon \, \phi^{\prime}_{i}(b) \Bigr) u_b(t) 
\end{align*}
for $i = 1, 2, 3$. That is, a nonhomogeneous fully nonlinear elliptic 
problem is solved at each time step.

\section{Numerical experiments}\label{sec-5}

In this section, we present a series of numerical tests to demonstrate 
the utility of the proposed mIP-DG methods for fully nonlinear PDEs of
the types \eqref{pde_ell} and $(\ref{pde})$. In all of our tests we shall use 
uniform spatial meshes as well as uniform temporal meshes for the dynamic 
problems. To solve the resulting nonlinear algebraic systems, we use 
the Matlab built-in nonlinear solver {\em fsolve} for the job.
For the elliptic problems we choose the initial guess as the 
linear interpolant of the boundary data $u_a$ and $u_b$. 
For dynamic problems, we let $u^0_h = \cP_h u_0$,  
$p_{1h}^0 = u^{0}_{hxx} (x^-)$, $p_{2h}^0 = \{ u^{0}_{hxx} (x)\}$, 
and  $p_{3h}^0 = u^{0}_{hxx} (x^+)$. Also, the initial guess for $u^n_h$ 
will be provided by $u^{n-1}_h$, and the initial guesses for 
$p_{1h}^n$, $p_{2h}^n$, and $p_{3h}^n$ will be provided by 
$p_{1h}^{n-1}$, $p_{2h}^{n-1}$, and $p_{3h}^{n-1}$, respectively.
For convenience, we set $\epsilon = 0$ for all tests. We remark that 
similar results can be obtained when $\epsilon \neq 0$, and 
the actual benefit of the symmetrization parameter is unclear
in the context of nonlinear algebraic systems. The role of $\alpha$ 
and the numerical moment will be further explored in section~\ref{alpha_tests}.

For our numerical tests, errors will be measured in the $L^\infty$ norm 
and the $L^2$ norm, where the errors are measured at the current time 
step for the dynamic problems. For the dynamic test problems, 
we shall see that the lower order time discretization dominates
the approximation error for reasonable time step size $\Delta t$.  
For the elliptic test problems and for the dynamic test problems
where the error is dominated by the spatial discretizations, 
it appears that the spatial error is of order $\mathcal{O} (h^\ell)$, where
\[
\ell =  \begin{cases}
	r+1, & \text{for } r \text{ odd} , \\
	r, & \text{for } r \text{ even} .
	\end{cases} 
\]
Furthermore, we observe that when using odd order elements, the schemes
exhibit optimal rate of convergence in both norms.

\subsection{Elliptic test problems}
We first present the results for three test problems of type \eqref{pde_ell}.
Both Monge-Amp\`ere and Bellman types of equations will be tested. 

\smallskip
{\bf Test 1.} 
Consider the stationary Monge-Amp\`{e}re problem
\begin{align*}
- u_{x x}^2 + 1 & = 0 , \qquad 0 < x < 1 , \\
u(0) = 0, \quad  u(1) & = \frac12 .
\end{align*}
It is easy to check that this problem has exactly two classical solutions: 
\[
u^+ (x) = \frac{1}{2} x^2 , \qquad u^- (x) = - \frac{1}{2} x^2 + x , 
\]
where $u^+$ is convex and $u^-$ is concave.  Note that $u^+$ is the unique 
viscosity solution which we want our numerical schemes to converge to.
In section~\ref{alpha_tests} we shall give some insights about 
the selectiveness of our schemes. 

\begin{figure}
\centerline{
\includegraphics[scale=0.4]{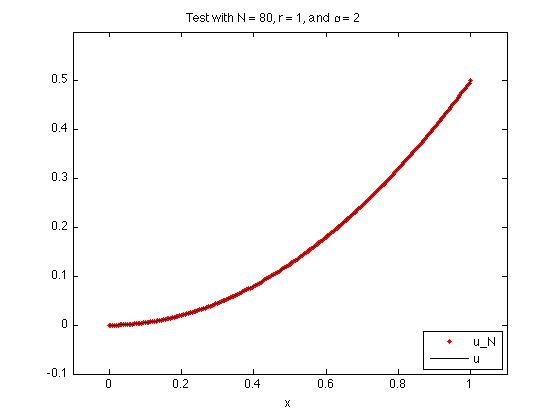}
}
\centerline{
\begin{tabular}{| c | c | c | c c | c c | c c |} 
		\hline
	$r$ & Norm & $h = 1/10$ 
		&\multicolumn{2}{| c |}{$h = 1/20$}
		&\multicolumn{2}{| c |}{$h = 1/40$}
		&\multicolumn{2}{| c |}{$h = 1/80$} \\ 
		\cline{3-9} 
	&   & Error & Error & Order & Error & Order & Error & Order  \\ 
		\hline \cline{1-9}
	1 & $L^{2}$ & 2.9e-03 & 7.3e-04 & 2.00 & 1.8e-04 
	  & 1.99 & 4.7e-05 & 1.97  \\ 
	  & $L^{\infty}$ & 3.8e-03 & 9.4e-04 & 2.00 & 2.4e-04 
	  & 1.99 & 6.1e-05 & 1.96  \\ 
		\hline
\end{tabular}
}
\caption{Test 1: $\epsilon = 0$, $\alpha = 2$, 
$\gamma_{01} = \gamma_{03} = 1$, and $\gamma_{02} = 1.1$.}
\label{exx_1_1_2}
\end{figure}

We approximate the given problem using the linear element ($r = 1$) 
to see how the approximation converges with respect to $h$ when the 
solution is not in the approximation space. The numerical results 
are shown in Figure~\ref{exx_1_1_2}. The results for the quadratic 
element ($r=2$) are presented in Figure~\ref{exx_1_1_1}. We note that
the approximations using $r=2$ are almost exact for each mesh size.  
This is expected since $u^+ \in V^h$ when $r = 2$.

\begin{figure} 
\centerline{
\includegraphics[scale=0.35]{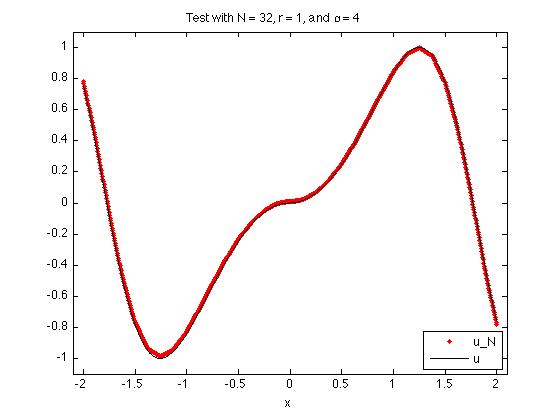}
\includegraphics[scale=0.35]{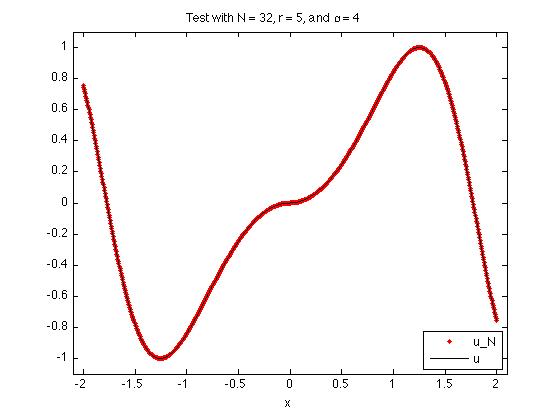}
}
\centerline{
\begin{tabular}{| c | c | c | c c | c c | c c |} 
		\hline
	$r$ & Norm & $h = 1$ 
		&\multicolumn{2}{| c |}{$h = 1/2$}
		&\multicolumn{2}{| c |}{$h = 1/4$}
		&\multicolumn{2}{| c |}{$h = 1/8$} \\ 
		\cline{3-9} 
	&   & Error & Error & Order & Error & Order & Error & Order  \\ 
		\hline \cline{1-9}
	1 & $L^{2}$ & 8.1e-01 & 2.4e-01 & 1.73 & 8.0e-02 
	  & 1.60 & 2.8e-02 & 1.52  \\ 
	  & $L^{\infty}$ & 1.0e+00 & 2.3e-01 & 2.14 & 7.8e-02 
	  & 1.58 & 2.7e-02 & 1.54  \\ 
		\hline
	2 & $L^{2}$ & 1.1e+00 & 2.9e-01 & 1.88 & 4.2e-02 
	  & 2.78 & 2.9e-02 & 0.56  \\ 
	  & $L^{\infty}$ & 8.1e-01 & 2.4e-01 & 1.76 & 4.5e-02 
	  & 2.40 & 1.8e-02 & 1.30  \\ 
		\hline
	3 & $L^{2}$ & 6.4e-01 & 2.7e-02 & 4.55 & 1.4e-03 
	  & 4.33 & 6.5e-05 & 4.38  \\ 
	  & $L^{\infty}$ & 4.9e-01 & 3.1e-02 & 3.99 & 1.6e-03 
	  & 4.32 & 9.1e-05 & 4.09  \\ 
		\hline
	4 & $L^{2}$ & 5.6e-02 & 3.2e-03 & 4.14 & 2.4e-04 
	  & 3.72 & 1.7e-05 & 3.83  \\ 
	  & $L^{\infty}$ & 4.9e-02 & 3.0e-03 & 4.02 & 2.6e-04 
	  & 3.56 & 1.6e-05 & 4.02  \\ 
		\hline
	5 & $L^{2}$ & 2.3e-02 & 8.5e-04 & 4.79 & 1.5e-05 
	  & 5.82 & 2.4e-07 & 5.96  \\ 
	  & $L^{\infty}$ & 2.1e-02 & 9.3e-04 & 4.49 & 1.8e-05 
	  & 5.67 & 2.6e-07 & 6.11  \\ 
		\hline
\end{tabular}
}
\caption{Test 2:  $\epsilon = 0$, $\alpha = 4$, 
$\gamma_{01} = \gamma_{03} = 2$, and $\gamma_{02} = 2.5$.}
\label{exx_1_2}
\end{figure}

\medskip
{\bf Test 2.} 
Consider the problem
\begin{align*}
- u_{ x x}^3 + \left| u_x \right| + S(x) & = 0, \qquad -2 < x < 2, \\ 
u(-2) = \sin (4),\quad u(2) & = - \sin (4) ,     
\end{align*}
where
\[
S(x)=\bigl[2 \mbox{sign}(x) \cos( x^2 ) - 4 \, x^2 \sin (x |x|) \bigr]^3 
       -2 |x\cos (x^2)| .
\]
This problem has the exact (viscosity) solution $u(x) =\sin (x|x|)$. Notice
that the equation is nonlinear in both $u_{x x}$ and $u_x$, and 
the exact solution is not twice differentiable at $x = 0$.  
The numerical results are shown in Figure~\ref{exx_1_2}.
As expected, we can see from the plot that the error appears largest 
around the point $x = 0$, and both the accuracy and order of 
convergence improve as the order of the element increases.  

\begin{figure} 
\centerline{
\includegraphics[scale=0.35]{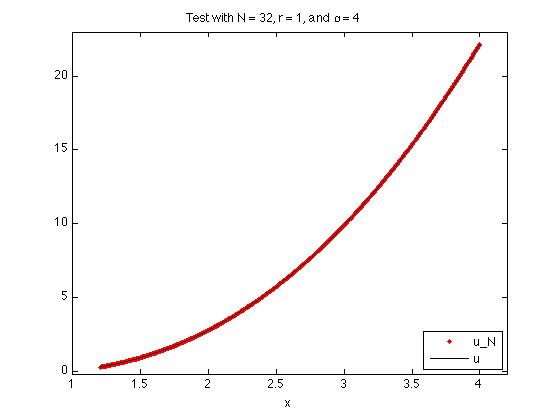}
\includegraphics[scale=0.35]{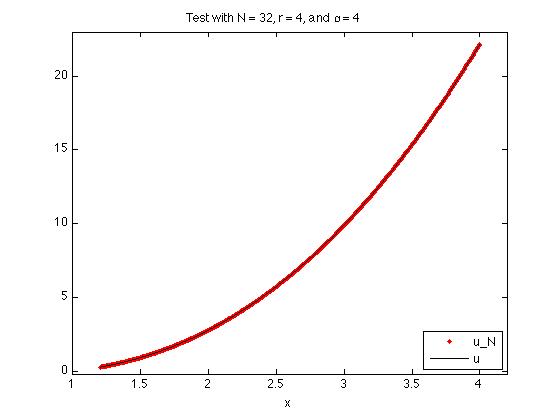}
}
\centerline{
\begin{tabular}{| c | c | c | c c | c c | c c |} 
		\hline
	$r$ & Norm & $h = 2.8/4$ 
		&\multicolumn{2}{| c |}{$h = 2.8/8$}
		&\multicolumn{2}{| c |}{$h = 2.8/16$}
		&\multicolumn{2}{| c |}{$h = 2.8/32$} \\ 
		\cline{3-9} 
	&   & Error & Error & Order & Error & Order & Error & Order  \\ 
		\hline \cline{1-9}
	1 & $L^{2}$ & 3.5e-01 & 9.8e-02 & 1.83 & 2.6e-02 
	  & 1.93 & 6.6e-03 & 1.97  \\ 
	  & $L^{\infty}$ & 3.9e-01 & 1.2e-01 & 1.70 & 3.4e-02 
	  & 1.81 & 9.0e-03 & 1.91  \\ 
		\hline
	2 & $L^{2}$ & 9.1e-03 & 1.9e-03 & 2.28 & 4.2e-04 
	  & 2.18 & 9.6e-05 & 2.11  \\ 
	  & $L^{\infty}$ & 9.9e-03 & 1.7e-03 & 2.53 & 3.6e-04 
	  & 2.23 & 8.2e-05 & 2.15  \\ 
		\hline
	3 & $L^{2}$ & 3.5e-04 & 2.7e-05 & 3.69 & 1.9e-06 
	  & 3.85 & 4.2e-07 & 2.14  \\ 
	  & $L^{\infty}$ & 5.1e-04 & 4.2e-05 & 3.61 & 3.3e-06 
	  & 3.69 & 3.7e-07 & 3.15  \\ 
		\hline
	4 & $L^{2}$ & 2.5e-05 & 1.4e-06 & 4.14 & 7.7e-08 
	  & 4.19 & 8.5e-09 & 3.18  \\ 
	  & $L^{\infty}$ & 3.3e-05 & 1.5e-06 & 4.46 & 7.6e-08 
	  & 4.30 & 1.3e-08 & 2.51  \\ 
		\hline
\end{tabular}
}
\caption{Test 3: $\epsilon = 0$, $\alpha = 4$, 
$\gamma_{01} = \gamma_{03} = 2$, and $\gamma_{02}= 2.5$.}
\label{exx_1_3}
\end{figure}

\medskip
{\bf Test 3.} 
Consider the stationary Hamilton-Jacobi-Bellman problem 
\begin{align*}
\inf_{0 < \theta(x) \leq 1} \left\{ 
- \theta u_{ x x } + \theta^2 \, x^2 \, u_{x} + \frac{1}{x} u + S(x) \right\}
& = 0 , \qquad 1.2 < x < 4 , \\
u(1.2) = 1.44 \ln 1.2, \quad  u(4) & = 16 \ln 4 ,
\end{align*}
where
\[
S(x) = \frac{4 \ln(x)^2 + 12 \ln(x) + 9 - 8 x^4 \ln(x)^2 
- 4 x^4 \ln(x)}{4 x^3 \left[ 2 \ln(x) + 1 \right]} .
\]
It can be shown that the exact (viscosity) solution of this problem is 
given by $u(x) = x^2 \ln x$, which occurs when 
$\theta^* (x) = \frac{2 \ln(x)+3}{2x^3 [2\ln(x)+ 1]}$.
We solve this problem using various order elements and record 
the numerical results in Figure~\ref{exx_1_3}, which shows
that our mIP-DG methods also can handle the Bellman-type fully 
nonlinear PDEs very well.

\subsection{Parabolic test problems}
We now implement the proposed fully discrete forward and backward Euler
mIP-DG methods for approximating fully nonlinear parabolic equations
of the form \eqref{pde}. While the above formulation makes no attempt
to formally quantify a CFL condition for the forward Euler method, 
our test problems generally require $\Delta t= \mathcal{O} (h^2)$ 
to ensure the stability. In fact, the constant for the CFL condition 
appears to decrease as the order of the element increases.
Below we implement both the implicit and explicit methods for each test 
problem. However, we make no attempt to classify and compare the 
efficiency of the two methods.  Instead, we focus on testing and demonstrating
the usability of both fully discrete schemes and 
their promising potentials. For explicit tests, we record the 
parameter $\kappa_t$ which serves as the scaling constant
for the CFL condition, so we have $\Delta t = \kappa_t h^2$.
For implicit tests, we record computed solution with
various time step $\Delta t$.

\begin{figure} 
\centerline{
\includegraphics[scale=0.35]{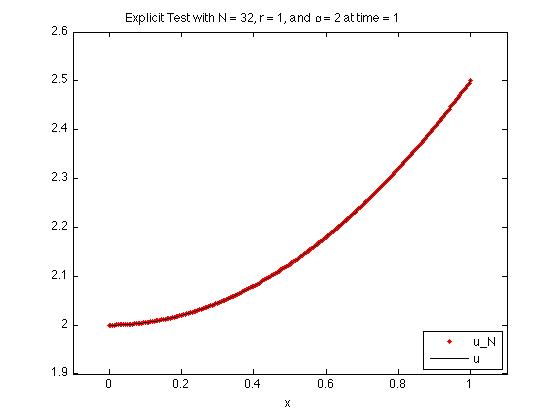}
\includegraphics[scale=0.35]{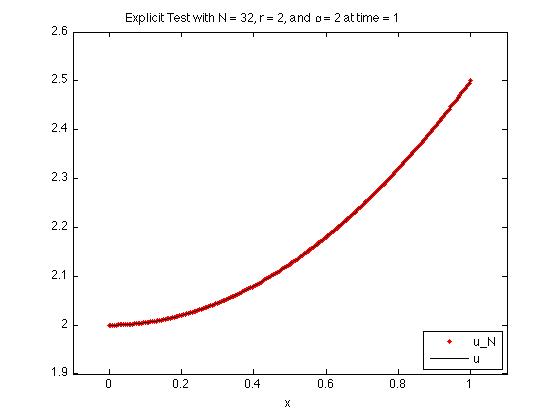}
}
\centerline{
\begin{tabular}{| c | c | c | c c | c c | c c |} 
		\hline
	$r$ & Norm & $h = 1/4$ 
		&\multicolumn{2}{| c |}{$h = 1/8$}
		&\multicolumn{2}{| c |}{$h = 1/16$}
		&\multicolumn{2}{| c |}{$h = 1/32$} \\ 
		\cline{3-9} 
	&   & Error & Error & Order & Error & Order & Error & Order  \\ 
		\hline \cline{1-9}
	1 & $L^{2}$ & 5.7e-03 & 1.4e-03 & 1.98 & 3.7e-04 
	  & 1.99 & 9.2e-05 & 1.99  \\ 
	  & $L^{\infty}$ & 7.9e-03 & 2.0e-03 & 1.99 & 5.0e-04 
	  & 1.99 & 1.3e-04 & 1.99  \\ 
		\hline
	2 & $L^{2}$ & 3.3e-05 & 8.2e-06 & 2.00 & 2.1e-06 
	  & 2.00 & 5.1e-07 & 2.00  \\ 
	  & $L^{\infty}$ & 4.5e-05 & 1.1e-05 & 2.00 & 2.8e-06 
	  & 2.00 & 7.1e-07 & 2.00  \\ 
		\hline
	3 & $L^{2}$ & 3.3e-05 & 8.2e-06 & 2.00 & 2.1e-06 
	  & 2.00 & 5.1e-07 & 2.00  \\ 
	  & $L^{\infty}$ & 4.5e-05 & 1.1e-05 & 2.00 & 2.8e-06 
	  & 2.00 & 7.1e-07 & 2.00  \\ 
		\hline
\end{tabular}
}
\caption{Test 4: Computed solutions at $T = 1$. 
$\kappa_t = 0.002$, $\epsilon = 0$, $\alpha = 2$, 
$\gamma_{01} = \gamma_{03} = 2$, and $\gamma_{02} = 2.5$.
Note, the scheme is unstable for $r = 2,3$ when $\kappa_t = 0.01$.}
\label{exx_2_1_e}
\end{figure}

\begin{figure} 
\centerline{
\includegraphics[scale=0.35]{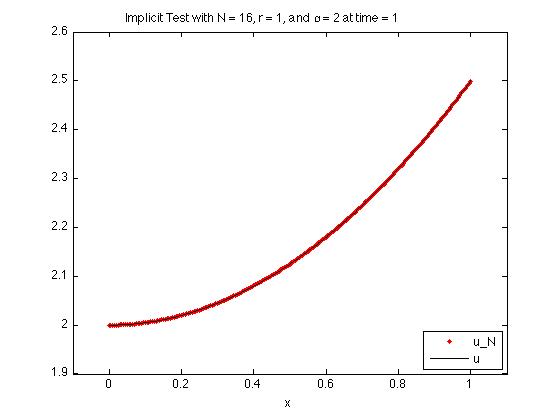}
\includegraphics[scale=0.35]{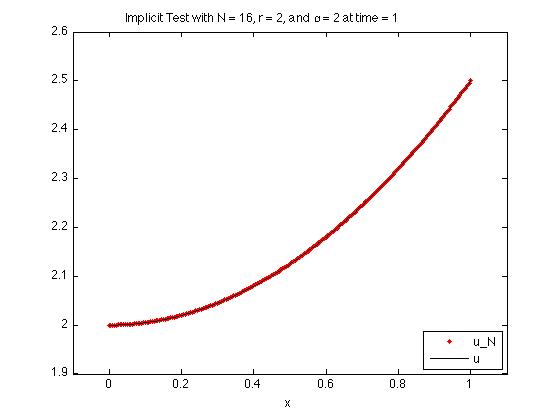} 
}
\centerline{
\begin{tabular}{| c | c | c | c c | c c |} 
		\hline
	$r$ & Norm & $h = 1/4$ 
		&\multicolumn{2}{| c |}{$h = 1/8$}
		&\multicolumn{2}{| c |}{$h = 1/16$} \\ 
		\cline{3-7} 
	&   & Error & Error & Order & Error & Order  \\ 
		\hline \cline{1-7}
	1 & $L^{2}$ & 4.4e-03 & 9.6e-04 & 2.20 & 1.8e-04 
	  & 2.40  \\ 
	  & $L^{\infty}$ & 9.4e-03 & 2.4e-03 & 2.00 & 5.9e-04 
	  & 2.00  \\ 
		\hline
	2 & $L^{2}$ & 2.6e-04 & 2.6e-04 & -0.00 & 2.6e-04 
	  & -0.00  \\ 
	  & $L^{\infty}$ & 3.6e-04 & 3.6e-04 & -0.00 & 3.6e-04 
	  & -0.00  \\ 
		\hline
	3 & $L^{2}$ & 2.6e-04 & 2.6e-04 & -0.00 & 2.6e-04 
	  & -0.00  \\ 
	  & $L^{\infty}$ & 3.6e-04 & 3.6e-04 & -0.00 & 3.6e-04 
	  & -0.00  \\ 
		\hline
\end{tabular}
}
\caption{Test 4: Computed solution at $T = 1$. 
$\Delta t = 0.001$, $\epsilon = 0$, $\alpha = 2$, 
$\gamma_{01} = \gamma_{03} = 2$, and $\gamma_{02} = 2.5$.}
\label{exx_2_1_i}
\end{figure}

\medskip
{\bf Test 4.} 
Let $\Omega = (0,1)$, $u_a(t) = t^4, u_b = \frac{1}{2} + t^4$, 
and $u_0(x) = \frac{1}{2} x^2$.  We consider the problem 
\eqref{pde}, \eqref{ic}, and \eqref{bc} with
\[
F(u_{xx}, u_x, u, t, x) = -u_{x x} \, u + \frac{1}{2} x^2 + t^4 - 4 \, t^3 + 1.
\]
It is easy to verify that this problem has a unique classical solution
$u(x,t) = 0.5 \, x^2 + t^4 + 1$. Notice that the PDE has a product 
nonlinearity in the second order derivative. The numerical results
for the fully discrete forward Euler method are presented in 
Figure~\ref{exx_2_1_e} and the results for the backward
Euler method are shown in Figure~\ref{exx_2_1_i}. We observe
that the errors for the backward Euler method are dominated by the 
relatively small size of the time step when compared to the forward 
Euler method. For smaller time step sizes, the errors are
similar. However, the backward Euler method appears unstable 
for $\kappa_t > 0.01$.

\begin{figure}
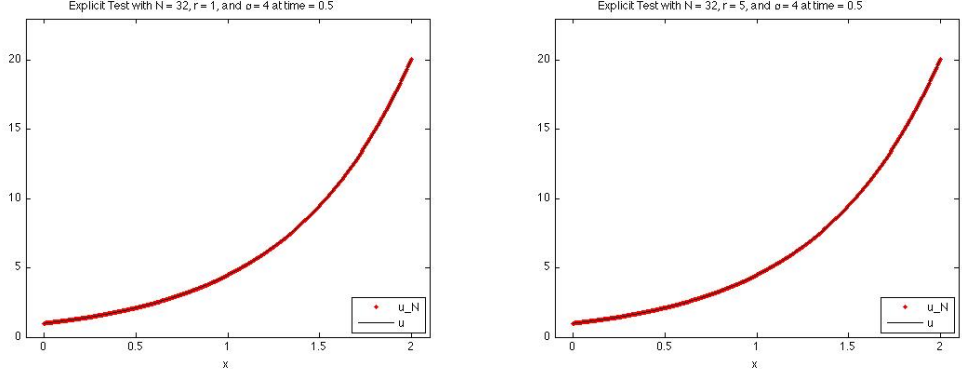
 
\centerline{
\begin{tabular}{| c | c | c | c c | c c | c c |} 
		\hline
	$r$ & Norm & $\Delta t = 1/10$ 
		&\multicolumn{2}{| c |}{$\Delta t = 1/20$}
		&\multicolumn{2}{| c |}{$\Delta t = 1/40$}
		&\multicolumn{2}{| c |}{$\Delta t = 1/80$} \\ 
		\cline{3-9} 
	&   & Error & Error & Order & Error & Order & Error & Order  \\ 
		\hline \cline{1-9}
	2 & $L^{2}$ & 2.4e-02 & 1.3e-02 & 0.93 & 6.4e-03 
	  & 0.96 & 3.2e-03 & 0.98  \\ 
	  & $L^{\infty}$ & 3.3e-02 & 1.7e-02 & 0.93 & 8.8e-03 
	  & 0.96 & 4.5e-03 & 0.98  \\ 
		\hline
\end{tabular}
}
\caption{Test 4: Computed solutions at $T = 1$.  
$h = 1/4$, $\epsilon = 0$, $\alpha = 2$, 
$\gamma_{01} = \gamma_{03}= 1$, and $\gamma_{02} = 1.1$.}
\label{exx_2_1_time_i}
\end{figure}

We now consider the error for the approximation resulting from 
using Euler time stepping methods. Note that the solution $u$ is a 
quadratic in space.  Letting $r = 2$, we limit the approximation 
error almost entirely to the time discretization scheme.  In fact, 
setting $t = 0$ and solving the stationary form of the PDE, we have 
\[
\| u - u_h\|_{L^2 ((0, 1))} \approx 1.6 \times 10^{-9} \quad \text{and} \quad
\| u - u_h\|_{L^\infty((0, 1))} \approx 2.4 \times 10^{-9}
\]
using the elliptic solver with $h = 1/4$, $\alpha = 2$, $\gamma_{01} = 
\gamma_{03} = 1$,  
$\gamma_{02}=1.1$, and initial guess given by the secant line for the boundary data.
Then, approximating the problem for varying $\Delta t$, we have the results recorded
in Figure~\ref{exx_2_1_time_e} for the forward Euler method and in
Figure~\ref{exx_2_1_time_i} for the backward Euler method.  
We observe that the convergence rate in time appears to have order 1 as expected.

\begin{figure} 
\centerline{
\begin{tabular}{| c | c | c | c c | c c | c c |} 
		\hline
	$r$ & Norm & $\kappa_t = 0.008$ 
		&\multicolumn{2}{| c |}{$\kappa_t = 0.004$}
		&\multicolumn{2}{| c |}{$\kappa_t = 0.002$}
		&\multicolumn{2}{| c |}{$\kappa_t = 0.001$} \\ 
		\cline{3-9} 
	&   & Error & Error & Order & Error & Order & Error & Order  \\ 
		\hline \cline{1-9}
	2 & $L^{2}$ & 8.2e-06 & 4.1e-06 & 1.00 & 2.1e-06 
	  & 1.00 & 1.0e-06 & 1.00  \\ 
	  & $L^{\infty}$ & 1.1e-05 & 5.7e-06 & 1.00 & 2.8e-06 
	  & 1.00 & 1.4e-06 & 1.00  \\ 
		\hline
\end{tabular}
}
\caption{Test 4: Computed solutions at time $T = 1$.
$h = 1/16$, $\epsilon = 0$, $\alpha = 2$, 
$\gamma_{01} = \gamma_{03} = 1$, and $\gamma_{02} = 1.1$.}
\label{exx_2_1_time_e}
\end{figure}

\begin{figure} 
\centerline{
\includegraphics[scale=0.35]{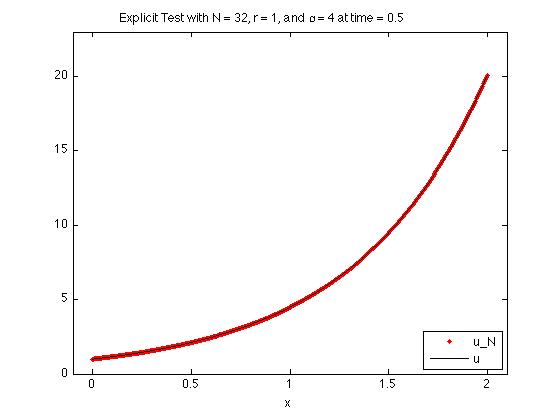} 
\includegraphics[scale=0.35]{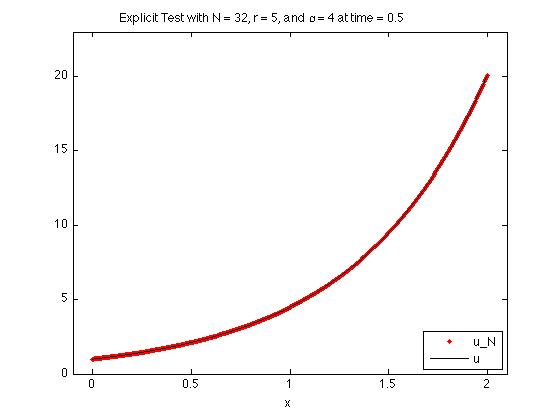}
}
\centerline{
\begin{tabular}{| c | c | c | c c | c c | c c |} 
		\hline
	$r$ & Norm & $h = 1/2$ 
		&\multicolumn{2}{| c |}{$h = 1/4$}
		&\multicolumn{2}{| c |}{$h = 1/8$}
		&\multicolumn{2}{| c |}{$h = 1/16$} \\ 
		\cline{3-9} 
	&   & Error & Error & Order & Error & Order & Error & Order  \\ 
		\hline \cline{1-9}
	1 & $L^{2}$ & 5.0e-01 & 3.6e-02 & 1.73 & 1.2e-02 
	  & 1.32 & 3.6e-03 & 1.67  \\ 
	  & $L^{\infty}$ & 8.2e-01 & 2.8e-01 & 1.57 & 1.0e-01 
	  & 1.47 & 3.1e-02 & 1.69  \\ 
		\hline
	2 & $L^{2}$ & 4.5e-02 & 1.2e-02 & 1.89 & 3.3e-03 
	  & 1.87 & 8.7e-04 & 1.93  \\ 
	  & $L^{\infty}$ & 6.0e-02 & 1.4e-02 & 2.11 & 3.6e-03 
	  & 1.96 & 9.0e-04 & 1.98  \\ 
		\hline
	3 & $L^{2}$ & 1.5e-03 & 2.8e-04 & 2.39 & 7.1e-05 
	  & 1.98 & 1.8e-05 & 1.98  \\ 
	  & $L^{\infty}$ & 2.7e-03 & 3.5e-04 & 2.97 & 7.6e-05 
	  & 2.21 & 1.8e-05 & 2.05  \\ 
		\hline
	4 & $L^{2}$ & 1.2e-03 & 2.9e-04 & 2.06 & 7.2e-05 
	  & 2.02 & 1.8e-05 & 2.01  \\ 
	  & $L^{\infty}$ & 1.3e-03 & 3.0e-04 & 2.13 & 7.3e-05 
	  & 2.02 & 1.8e-05 & 2.01  \\ 
		\hline
	5 & $L^{2}$ & 1.2e-03 & 2.9e-04 & 2.00 & 7.2e-05 
	  & 2.00 & 1.8e-05 & 2.00  \\ 
	  & $L^{\infty}$ & 1.2e-03 & 2.9e-04 & 2.00 & 7.3e-05 
	  & 2.00 & 1.8e-05 & 2.00  \\ 
		\hline
\end{tabular}
}
\caption{Test 5:  Computed solutions at time $T = 3.10$.  
$\kappa_t = 0.0025$, $\epsilon = 0$, $\alpha = 4$, 
$\gamma_{01} = \gamma_{03} = 2$, and $\gamma_{02} = 2.5$.
The method is not stable for $\kappa_t = 0.005$.}
\label{exx_2_2_e}
\end{figure}

\begin{figure} 
\centerline{
\includegraphics[scale=0.35]{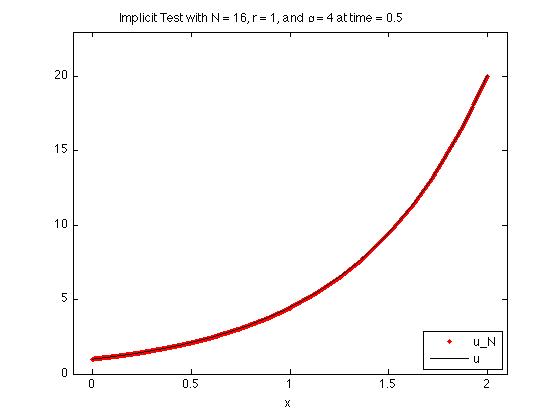}
\includegraphics[scale=0.35]{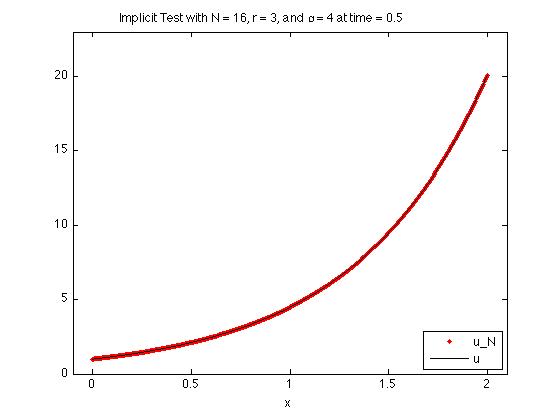}
}
\centerline{
\begin{tabular}{| c | c | c | c c | c c |} 
		\hline
	$r$ & Norm & $h = 1/2$ 
		&\multicolumn{2}{| c |}{$h = 1/4$}
		&\multicolumn{2}{| c |}{$h = 1/8$} \\ 
		\cline{3-7} 
	&   & Error & Error & Order & Error & Order \\ 
		\hline \cline{1-7}
	1 & $L^{2}$ & 4.2e-01 & 1.4e-01 & 1.54 & 4.6e-02 
	  & 1.66 \\ 
	  & $L^{\infty}$ & 8.3e-01 & 2.4e-01 & 1.77 & 7.9e-02 
	  & 1.63 \\ 
		\hline
	2 & $L^{2}$ & 7.3e-02 & 1.6e-02 & 2.21 & 3.0e-03 
	  & 2.40 \\ 
	  & $L^{\infty}$ & 9.6e-02 & 1.8e-02 & 2.41 & 3.2e-03 
	  & 2.49 \\ 
		\hline
	3 & $L^{2}$ & 2.8e-03 & 7.8e-04 & 1.82 & 9.1e-04 
	  & -0.22 \\ 
	  & $L^{\infty}$ & 5.6e-03 & 8.5e-04 & 2.71 & 9.2e-04 
	  & -0.11 \\ 
		\hline
\end{tabular}
}
\caption{Test 5: Computed solutions at time $T = 0.5$. 
$\Delta t = 0.0005$, $\epsilon = 0$, $\alpha = 4$, 
$\gamma_{01} = \gamma_{03} = 2$, and $\gamma_{02} = 2.5$.}
\label{exx_2_2_i}
\end{figure}

\medskip
{\bf Test 5.} 
Let $\Omega = (0,2)$, $u_a(t) = 1$, $u_b = e^{2 (t+1)}$, and $u_0(x) = e^x$.  
We consider the problem \eqref{pde}, \eqref{ic}, and \eqref{bc} with
\[
F(u_{x x}, u_x, u, t, x) = - u_x \ln \bigl( u_{x x} + 1 \bigr) + S(x,t) , 
\]
and
\[
S(x,t) = e^{(t+1) x} \Big(x - (t+1) \ln \bigl((t+1)^2 e^{(t+1) x}+1 \bigr) \Big).
\]
It is easy to verify that this problem has a unique classical
solution $u(x,t) = e^{(t+1) x}$. Notice that this problem is nonlinear in 
both $u_{x x}$ and $u_x$. Furthermore, the exact solution $u$ cannot be 
factored into the form $u(x,t) = G(t) \, Y (x)$ for some functions $G$ and $Y$.  
The numerical results for the fully discrete forward Euler method
are recorded in Figure~\ref{exx_2_2_e} and the results for the backward
Euler method are given in Figure~\ref{exx_2_2_i}. The error appears to be 
dominated by the low order time discretization given the relatively
large value for $\Delta t$ in the backward Euler test. However, using 
a smaller $\Delta t$ for the forward Euler test, we were able to achieve 
a higher order of accuracy.  We remark that even for $\Delta t = 0.005 h^2$, 
the forward Euler scheme is not stable for $h = \frac14$ and $r = 1$.

\begin{figure} 
\centerline{
\includegraphics[scale=0.35]{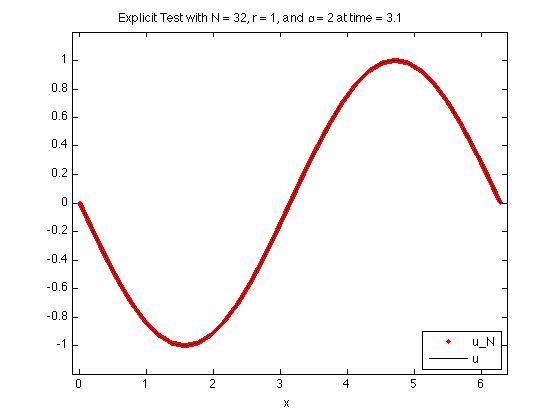}
\includegraphics[scale=0.35]{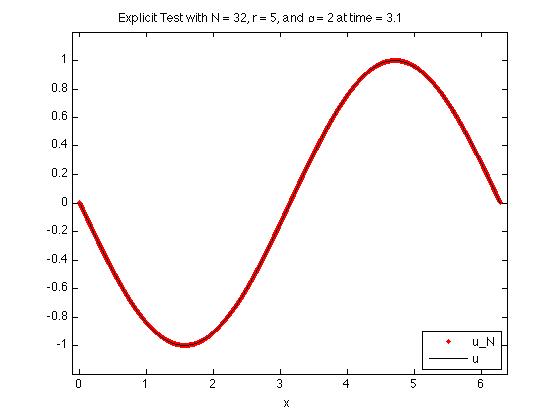}
}
\centerline{
\begin{tabular}{| c | c | c | c c | c c | c c |} 
		\hline
	$r$ & Norm & $h = \pi/2$ 
		&\multicolumn{2}{| c |}{$h = \pi/4$}
		&\multicolumn{2}{| c |}{$h = \pi/8$}
		&\multicolumn{2}{| c |}{$h = \pi/16$} \\ 
		\cline{3-9} 
	&   & Error & Error & Order & Error & Order & Error & Order  \\ 
		\hline \cline{1-9}
	1 & $L^{2}$ & 2.2e-01 & 5.3e-02 & 2.07 & 1.3e-02 
	  & 2.02 & 3.3e-03 & 2.01  \\ 
	  & $L^{\infty}$ & 1.7e-01 & 4.8e-02 & 1.87 & 1.2e-02 
	  & 1.98 & 3.1e-03 & 1.99  \\ 
		\hline
	2 & $L^{2}$ & 6.0e-02 & 1.6e-02 & 1.90 & 4.2e-03 
	  & 1.94 & 1.1e-03 & 1.97  \\ 
	  & $L^{\infty}$ & 6.4e-02 & 1.5e-02 & 2.07 & 3.5e-03 
	  & 2.13 & 8.2e-04 & 2.09  \\ 
		\hline
	3 & $L^{2}$ & 7.4e-03 & 6.9e-04 & 3.43 & 1.4e-04 
	  & 2.32 & 3.5e-05 & 2.00  \\ 
	  & $L^{\infty}$ & 8.0e-03 & 5.6e-04 & 3.82 & 1.0e-04 
	  & 2.46 & 2.3e-05 & 2.14  \\ 
		\hline
	4 & $L^{2}$ & 2.5e-03 & 5.7e-04 & 2.10 & 1.4e-04 
	  & 2.03 & 3.5e-05 & 2.01  \\ 
	  & $L^{\infty}$ & 1.4e-03 & 3.5e-04 & 2.01 & 8.9e-05 
	  & 1.98 & 2.2e-05 & 1.99  \\ 
		\hline
	5 & $L^{2}$ & 2.2e-03 & 5.6e-04 & 2.00 & 1.4e-04 
	  & 2.00 & 3.5e-05 & 2.00  \\ 
	  & $L^{\infty}$ & 1.4e-03 & 3.6e-04 & 1.99 & 8.9e-05 
	  & 2.00 & 2.2e-05 & 2.00  \\ 
		\hline
\end{tabular}
}
\caption{Test 6: Computed solutions at time $T = 3.10$. 
$\kappa_t = 0.002$, $\epsilon = 0$, $\alpha = 2$, 
$\gamma_{01} = \gamma_{03} = 2$, and $\gamma_{02} = 2.5$.}
\label{exx_2_3_e}
\end{figure}

\begin{figure}
\centerline{
\includegraphics[scale=0.35]{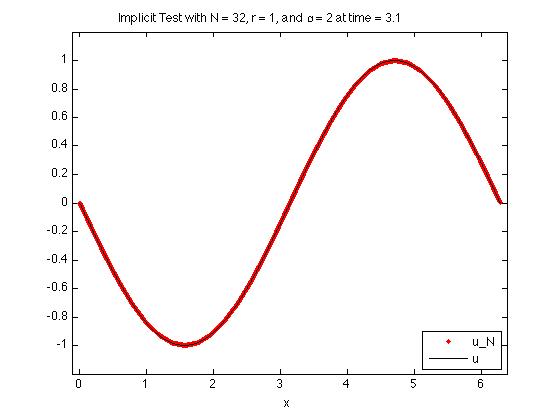}
\includegraphics[scale=0.35]{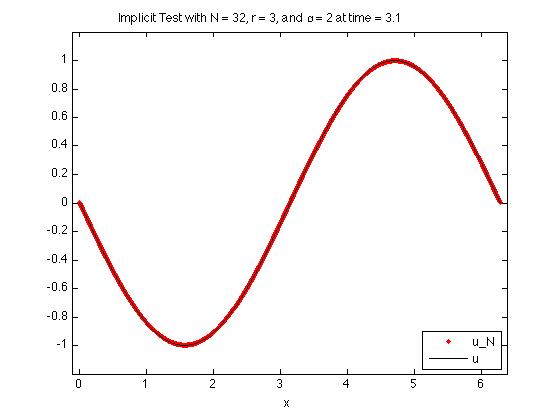}
}
\centerline{
\begin{tabular}{| c | c | c | c c | c c | c c |} 
		\hline
	$r$ & Norm & $h = \pi/2$ 
		&\multicolumn{2}{| c |}{$h = \pi/4$}
		&\multicolumn{2}{| c |}{$h = \pi/8$}
		&\multicolumn{2}{| c |}{$h = \pi/16$} \\ 
		\cline{3-9} 
	&   & Error & Error & Order & Error & Order & Error & Order  \\ 
		\hline \cline{1-9}
	1 & $L^{2}$ & 1.7e-01 & 4.9e-02 & 1.82 & 1.4e-02 
	  & 1.84 & 4.8e-03 & 1.50  \\ 
	  & $L^{\infty}$ & 1.5e-01 & 4.4e-02 & 1.78 & 1.3e-02 
	  & 1.82 & 4.1e-03 & 1.60  \\ 
		\hline
	2 & $L^{2}$ & 8.0e-02 & 2.0e-02 & 2.00 & 5.9e-03 
	  & 1.76 & 3.2e-03 & 0.87  \\ 
	  & $L^{\infty}$ & 7.0e-02 & 1.6e-02 & 2.14 & 4.0e-03 
	  & 1.98 & 1.9e-03 & 1.06  \\ 
		\hline
	3 & $L^{2}$ & 1.1e-02 & 3.0e-03 & 1.91 & 2.8e-03 
	  & 0.09 & 2.8e-03 & 0.00  \\ 
	  & $L^{\infty}$ & 8.1e-03 & 1.8e-03 & 2.16 & 1.8e-03 
	  & 0.01 & 1.8e-03 & 0.00  \\ 
		\hline
\end{tabular}
}
\caption{Test 6: Computed solutions at time $T = 3.10$. 
$\Delta t = 0.0062$, $\epsilon = 0$, $\alpha = 2$, 
$\gamma_{01} = \gamma_{03} = 2$, and $\gamma_{02} = 2.5$.}
\label{exx_2_3_i}
\end{figure}

\medskip
{\bf Test 6.} 
Let $\Omega = (0, 2 \pi)$, $u_a(t) = 0$, $u_b = 0$, and $u_0(x) = \sin (x)$.  
We consider the problem \eqref{pde}, \eqref{ic}, and \eqref{bc} with
\[
F(u_{x x}, u_x, u, t, x) = - \min_{\theta (t, x) \in \{1,2 \}} 
\Big\{A_\theta\, u_{x x}-c(x,t) \cos (t) \, \sin (x) - \sin(t) \, \sin(x) \Big\} , 
\]
where $A_1 = 1, A_2 = \frac12$, and
\[
c(x,t) = \begin{cases}
1, & \text{if } 
0 < t \leq \frac{\pi}2 \mbox{ and } 0 < x \leq \pi \text{ or } 
\frac{\pi}2 < t \leq \pi \text{ and } \pi < x < 2 \pi, \\
\frac12, & \mbox{otherwise}. 
\end{cases}
\]
It is easy to check that this problem has a unique classical solution
$u (x,t) = \cos(t) \sin(x)$.  Notice that this problem has a finite 
dimensional control parameter set, and the optimal control is given by
\[
\theta^* (t, x) = \begin{cases}
1, & \mbox{if } c(x,t) = 1 , \\
2 , & \mbox{if } c(x,t) = 2 .
\end{cases} 
\]
The numerical results are recorded in Figure~\ref{exx_2_3_e} for the fully
discrete forward Euler method and in Figure~\ref{exx_2_3_i} for the 
backward Euler method. We observe that the accuracy of the implicit method 
appears to suffer from the lower order accuracy of the Euler method.  For 
$h = \frac{\pi}{8}$, the explicit method requires 
$\Delta t \approx 3.1 \times 10^{-4}$, while the implicit method 
only needs $\Delta t = 0.062$. When $\Delta t$ increases, the explicit 
method demonstrates instability.

\subsection{The role of the numerical moment} \label{alpha_tests}

We now discuss the role and utility of the numerical moment 
in forming an appropriate numerical operator. Consider the stationary 
Monge-Amp\`{e}re problem from Test 1, which has the following two solutions: 
\[
u^+ (x) = \frac{1}{2} x^2 , \qquad u^- (x) = - \frac{1}{2} x^2 + x , 
\]
where $u^+$ is convex and $u^-$ is concave. The solution $u^+$ is the unique viscosity 
solution. 

To demonstrate the role of the numerical moment, we approximate 
the given problem using $\alpha > 0$, $\alpha = 0$, and $\alpha < 0$.  
Notice that multiplying the PDE by $-1$, we see that $u^-$ is the unique 
viscosity solution of the equation with the operator $F(u) = u_{x x}^2 - 1$.
Then, for $\alpha > 0$, our scheme should converge to $u^+$, and 
for $\alpha < 0$ our scheme should converge to $u^-$ provided 
$|\alpha |$ is sufficiently large. However, for $\alpha = 0$, the scheme 
may converge to either $u^+$ or $u^-$ depending on the initial guess 
used for the nonlinear solver. Note that while we cannot globally 
bound $\partial_{u^{\prime \prime}} F$ for the operator 
$F(u^{\prime \prime})=1 -(u^{\prime \prime})^2$, we can
locally bound $\partial_{u^{\prime \prime}} F$. Thus, the necessary 
magnitude for $\alpha$ to allow selective convergence depends on the 
initial guess and the solver. Without a global bound on 
$\partial_{u^{\prime \prime}} F$, the numerical operator is only
locally monotone.

\begin{figure} 
\centerline{
\includegraphics[scale=0.35]{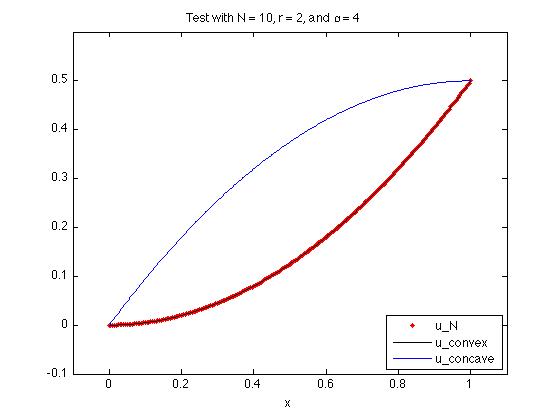}
\includegraphics[scale=0.35]{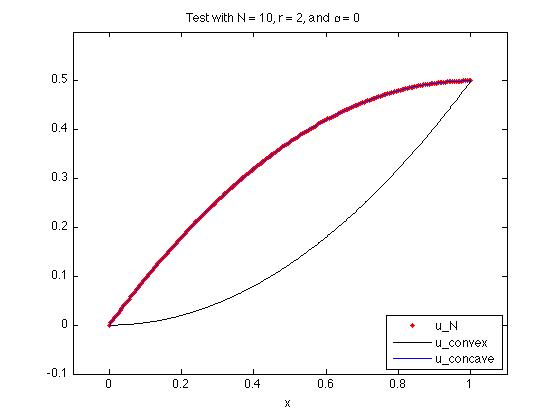}
}
\centerline{
\includegraphics[scale=0.4]{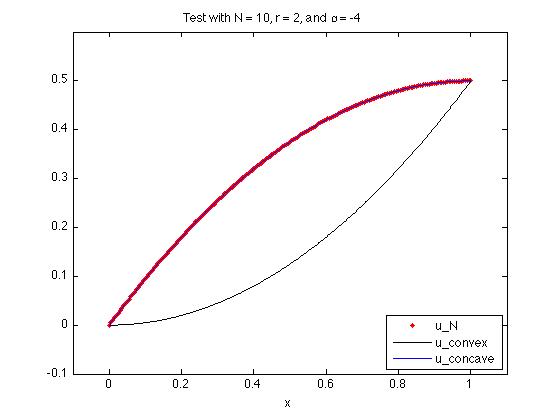}
}
\centerline{
\begin{tabular}{| c | c | c | c |}
		\hline
	Norm & $\alpha = 4$ & $\alpha = 0$ & $\alpha = -4$ \\
		\hline
	$L^{2}$ & 2.5e-08 & 5.3e-10 & 3.7e-10 \\
		\hline
	$L^\infty$ & 3.3e-08 & 8.6e-10 & 5.7e-10 \\
		\hline
\end{tabular}
}
\caption{Test 1:  $h = 1/10$, $\epsilon = 0$, $r = 2$, 
$\gamma_{01} = \gamma_{03} = 1.1$, and $\gamma_{02} = 1.5$.}
\label{exx_1_1_1}
\end{figure}

Let $\overline{u}$ be the linear interpolant of the boundary data 
and let the initial guess for $u$ be given by 
$u^{(0)} = \frac{1}{3} \overline{u} + \frac{2}{3} \, u^-$ and 
the initial guesses for $p_i$ be given by $p_i^{(0)} = 0$ 
for $i = 1, 2, 3$. Thus, the initial guess is closer to $u^-$.  
From Figure~\ref{exx_1_1_1} we see that the scheme converges 
to $u^+$ for $\alpha = 4$ and the scheme converges to $u^-$ 
for $\alpha = 0$ and $\alpha = -4$ for the given parameters.  
If we change the initial guess to 
$u^{(0)} = \frac{1}{3} \overline{u} + \frac{2}{3} \, u^+$, 
the scheme converges to $u^+$ for $\alpha = 0$ and $\alpha = 4$ 
and the scheme converges to $u^-$ for $\alpha = -4$ for the given parameters.  
Furthermore, for $u^{(0)} = \overline{u}$, {\em fsolve} does
not find a root for $\alpha = 0$, whereas the scheme converges 
for $\alpha = \pm 4$.

Therefore, the numerical moment plays two major roles.
It allows the scheme to converge for a wider range
of initial guesses, and it enables the scheme to address the 
issue of the conditional uniqueness of viscosity solutions.  
Given the form of the numerical moment, $\alpha (p_1 - 2 \, p_2 + p_3)$, 
these benefits are even more substantial given the way in 
which $p_1$, $p_2$, and $p_3$ are formed. The three variables only 
differ in their jump terms. When $\gamma_{01} = \gamma_{02} = \gamma_{03}$,
the three different choices for the numerical fluxes (or jump terms)
are all equivalent at the PDE level, and often the various 
jump formulations are presented as interchangeable when discretizing
linear and quasilinear PDEs using the DG methodology. Yet, for our schemes 
for fully nonlinear PDEs, we see that the three different choices 
of the numerical fluxes all play an essential role at the
numerical level when combined to form the numerical moment, even in 
the degenerative case where $\gamma_{01} = \gamma_{02} = \gamma_{03}$ 
which will be discussed below.

The role of the numerical moment can heuristically be understood as follows
when the numerical moment is rewritten in the form
\[
\alpha h^2 \Bigl(\frac{ p_1 - 2 p_2 + p_3}{h^2}\Bigr).
\]
From here, we can see that the numerical moment acts as a centered 
difference approximation for $u_{xxxx}$ multiplied by a factor 
that tends to zero with rate $\mathcal{O} \left( h^2 \right)$.
Thus, at the PDE level, we are in essence approximating 
the nonlinear elliptic operator
\[
F \left( u_{x x}, u_x, u, x \right)
\]
by the quasilinear fourth order operator $\widehat{F}_\rho$, where
\[
\widehat{F}_\rho \left( u_{x x x x}, u_{x x}, u_x, u, x \right) 
=\rho \, u_{x x x x} + F \left( u_{x x}, u_x, u, x \right) .
\]  
In the limit as $\rho \to 0$, we heuristically expect the unique limit
of the fourth order problem to converge to the unique viscosity solution 
of the second order problem.  Using a converging family of fourth order 
quasilinear PDEs to approximate a fully nonlinear second order PDE 
has previously been considered for PDEs such as the Monge-Amp\`{e}re equation, 
the prescribed Gauss curvature equation, the infinity-Laplace equation, 
and linear second order equations of non-divergence form.  The method is known
as the vanishing moment method. We refer the reader to 
\cite{Feng_Neilan11,FGN12} for a detailed exposition.

In addition to the connection with the numerical moment to quasilinear 
fourth order PDEs, we also mention another benefit of the numerical moment.  
By the choice of $\alpha$, we can enlarge the domain for which the 
numerical operator $\widehat{F}$ is increasing in $p_1$ and $p_3$ and 
decreasing in $p_2$.  Since the definition of ellipticity is based on 
the monotonicity of the operator, and the issue of conditional uniqueness
stems from whether the solution preserves the monotonicity of the 
operator, building monotonicity into the discretization is important 
when trying to preserve the nature of the operator we are approximating.  
	
We can demonstrate the power of the monotonicity of the numerical 
operator with two simple tests. For both tests, we shall again approximate 
the Monge-Amp\`{e}re problem from Test 1. However, now we let 
$\gamma_{01} = \gamma_{02} = \gamma_{03}$.  Then, we have 
$p_2 = \frac{p_1 + p_3}2$, which in turn implies that the 
equation for $p_2$ is redundant in the formulation and the numerical 
moment should be zero upon convergence to a root.

\begin{figure} 
\centerline{
\includegraphics[scale=0.35]{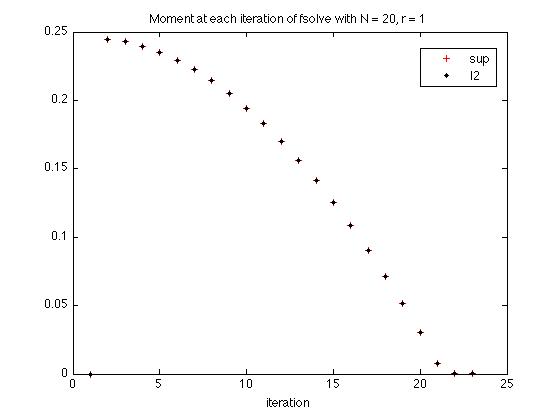}
\includegraphics[scale=0.35]{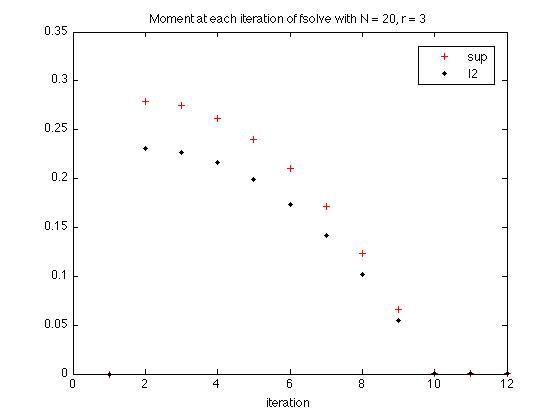}
}
\caption{Plots of the norm of $p_1 - 2 p_2 + p_3$ with $\alpha = 4$, 
$\epsilon = 0$, and $\gamma_{01} = \gamma_{02} = \gamma_{03} = 2$ 
at each iteration of {\em fsolve}.} \label{exx_1_1_moment}
\end{figure}

For the first test, we again approximate the Monge-Amp\`{e}re problem 
from Test 1 while plotting the norm of $p_1 - 2 p_2 + p_3$ after each 
iteration of {\em fsolve}. From Figure~\ref{exx_1_1_moment}, we can see that 
even though we expect the moment to be zero based on the redundancy of 
the equation for $p_2$ given the equations for $p_1$ and $p_3$, the Newton-based 
solver {\em fsolve} treats $p_1$, $p_2$, and $p_3$ as independent variables 
when searching for a root. The monotonicity of each variable appears 
to aid {\em fsolve} in the search for a root.

For the second test, instead of using {\em fsolve}, a Newton-based solver, 
for solving the nonlinear system of equations, we use the following 
splitting algorithm:

\medskip
\begin{solver} \label{solver_alg}
\noindent
\begin{enumerate}
\item[{\rm (1)}] Pick an initial guess for $u$, $p_1$, and $p_3$. 
\item[{\rm (2)}] Solve equation \eqref{pde_ell_weak} for $p_2$. 
\item[{\rm (3)}] Solve equation \eqref{pi_bil} for $i=2$ for $u$. 
\item[{\rm (4)}] Solve equation \eqref{pi_bil} for $i=1$ for $p_1$. 
\item[{\rm (5)}] Solve equation \eqref{pi_bil} for $i=3$ for $p_3$. 
\item[{\rm (6)}] Repeat Steps 2 - 5 until the change in $p_2$ is sufficiently small.
\end{enumerate}
\end{solver}

\smallskip
We observe that only Step (2) involves the use of a nonlinear solver.  
Each of Steps (3)-(5) only requires solving a linear system with a 
constant matrix that can be pre-factored. Thus, the above solver fully
decouples the entire system of equations and minimizes the number of 
unknowns in the nonlinear system. Because this paper is concerned mainly
with the discretization of fully nonlinear PDEs, we do not make
an effort to compare solvers. The simple solver presented here is meant 
to demonstrate a potential benefit of using the numerical moment to create 
monotone numerical operators.

\begin{figure}
\centerline{
\includegraphics[scale=0.35]{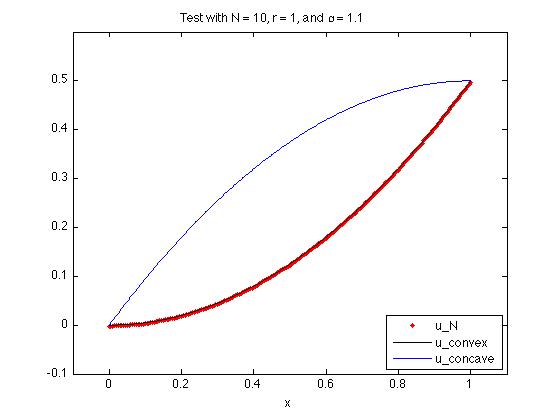}
\includegraphics[scale=0.35]{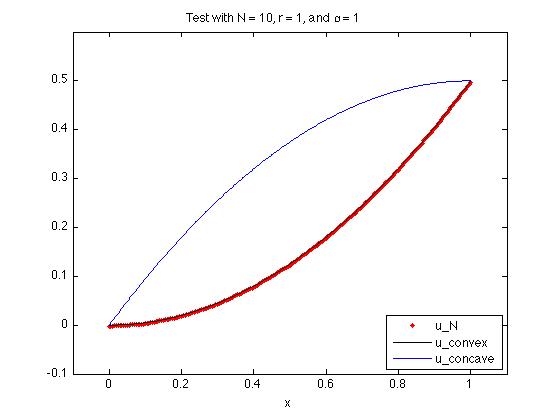}
}
\centerline{
\includegraphics[scale=0.35]{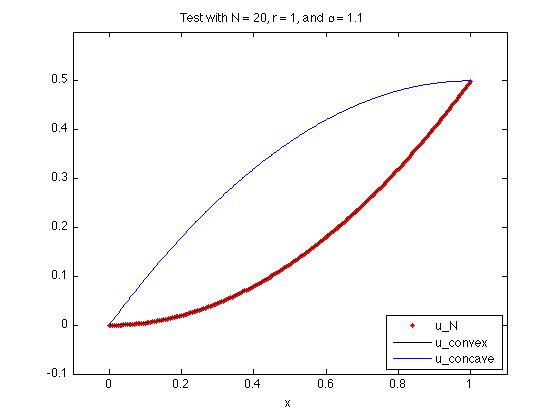}
\includegraphics[scale=0.35]{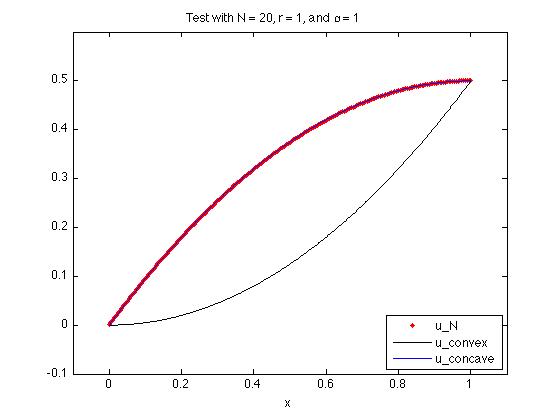}
}
\centerline{
\begin{tabular}{| c | c | c | c c | c c | c c |} 
		\hline
	$\alpha$ & Norm & $h = 1/5$ 
		&\multicolumn{2}{| c |}{$h = 1/10$}
		&\multicolumn{2}{| c |}{$h = 1/20$}
		&\multicolumn{2}{| c |}{$h = 1/40$} \\ 
		\cline{3-9} 
	&   & Error & Error & Order & Error & Order & Error & Order  \\ 
		\hline \cline{1-9}
	4 & $L^{2}$ & 6.8e-03 & 1.7e-03 & 2.00 & 4.3e-04 
	  & 2.00 & 1.0e-04 & 2.08  \\ 
	  & $L^{\infty}$ & 1.0e-02 & 2.5e-03 & 2.00 & 6.2e-04 
	  & 2.00 & 1.6e-04 & 2.00  \\ 
		\hline
	2 & $L^{2}$ & 6.8e-03 & 1.7e-03 & 2.00 & 4.3e-04 
	  & 2.01 & 9.8e-05 & 2.12  \\ 
	  & $L^{\infty}$ & 1.0e-02 & 2.5e-03 & 2.00 & 6.2e-04 
	  & 2.00 & 1.6e-04 & 2.00  \\ 
		\hline
	1.1 & $L^{2}$ & 6.8e-03 & 1.7e-03 & 2.00 & 4.2e-04 
	  & 2.01 & 9.7e-05 & 2.13  \\ 
	  & $L^{\infty}$ & 1.0e-02 & 2.5e-03 & 2.00 & 6.2e-04 
	  & 2.00 & 1.6e-04 & 2.00  \\ 
		\hline
	1 & $L^{2}$ & 6.8e-03 & 1.7e-03 & 2.00 & 5.7e-04 
	  & 1.58 & 8.2e-04 & -0.53  \\ 
	  & $L^{\infty}$ & 1.0e-02 & 2.5e-03 & 2.00 & 9.4e-04 
	  & 1.42 & 1.2e-03 & -0.32  \\ 
		\hline
	0.99 & $L^{2}$ & 6.0e-03 & 9.7e-04 & 2.62 & 5.7e-04 
	  & 0.77 & 8.2e-04 & -0.53  \\ 
	  & $L^{\infty}$ & 9.9e-03 & 2.5e-03 & 2.00 & 9.4e-04 
	  & 1.40 & 1.2e-03 & -0.32  \\ 
		\hline
	0 & $L^{2}$ & 6.8e-03 & 1.7e-03 & 2.00 & 4.3e-04 
	  & 1.99 & 1.1e-04 & 1.96  \\ 
	  & $L^{\infty}$ & 1.0e-02 & 2.5e-03 & 2.00 & 6.3e-04 
	  & 1.99 & 1.6e-04 & 1.96  \\ 
		\hline
\end{tabular}
}
\caption{Test 1 is solved using Algorithm~\ref{solver_alg} with 
$r = 1$, $\epsilon = 0$, and $\gamma_{01} = \gamma_{02} = \gamma_{03} = 2$.  
For $\alpha \geq 1.1$, the scheme converges to $u^+$.  
For $\alpha \leq 0.99$, the scheme converges to $u^-$.  
When $\alpha = 1.0$, the scheme converges to $u^+$ for $h \geq \frac{1}{10}$ and
the scheme converges to $u^-$ for $h \leq \frac{1}{20}$.}
\label{splitting}
\end{figure}

We use Algorithm~\ref{solver_alg} with {\em fsolve} to execute Step (2)
for Test 1. Let the initial guesses be given by $u = u^-$ and 
$p_1 = p_2 = p_3 = -0.99$. For $p_2 = -0.99$, $F$ is increasing 
while $\widehat{F}$ is decreasing for $\alpha > 0.99$. Since 
$F \left( -0.99 \right) > 0$ and $\widehat{F}$ is decreasing for 
$p_2 \geq -1$ when $\alpha > 1$, we expect the splitting algorithm 
will move away from the concave root $p_2= -1$. The numerical 
results are presented in Figure~\ref{splitting}. We note that even 
with the initial guess close to $u^-$, the solver, with the aid of 
the numerical moment, converges to $u^+$. Similarly, the solver 
converges to $u^+$ when $p_1 = p_2 = p_3 > -0.99$ are used as initial guesses.
For initial guesses $p_1 = p_2 = p_3 < -1.0$, the solver does not converge.  
Thus, we see that even for the above simple solver, the 
monotonicity of $\widehat{F}$ provided by the numerical moment allows 
the scheme to either selectively converge to $u^+$ or diverge and find 
no solution. Hence, we again see the benefit of including the numerical moment 
when tackling the issue of conditional uniqueness for viscosity solutions.

\section{Conclusion}\label{sec-6}
In this paper we present a general framework for constructing high order 
interior penalty discontinuous Galerkin methods for approximating 
viscosity solutions of fully nonlinear second order elliptic and parabolic PDEs.
The proposed framework extends the (second order) finite difference framework
developed by the authors in \cite{FKL11} to a more flexible DG framework, 
allowing approximating fully nonlinear PDEs using high order 
polynomials and non-uniform meshes. Various numerical experiments
are provided to show the performance of the proposed methodology.
The proposed DG framework is based on a nonstandard mixed formulation of 
the underlying fully nonlinear PDE. In order to capture discontinuities of 
the second order derivative $u_{xx}$ of the solution $u$, 
three independent functions $p_1, p_2$ and $p_3$ are 
introduced to accomplish the goal, where $p_1$ and $p_3$
measure the left and right limits of $u^{\prime \prime}$.
If $u^{\prime \prime}$ is discontinuous, $p_1$ and $p_3$ can be used to
gain insight into the discontinuity upon convergence. Thus,
the methodology has the ability to capture some of the more interesting
aspects of the viscosity solutions. The proposed mIP-DG methodology 
takes the most important aspects of the companion finite 
difference framework of \cite{FKL11} and extends 
them in multiple directions. For example, by adopting and expanding the idea 
of numerical operators, the mIP-DG formulation allows for even more 
flexibility than finite difference methods in construction.  

The proposed mIP-DG discretizations touch the inner core and make use of 
the full potential of the DG methodology. This is because there is a 
natural match among the three choices of possible numerical fluxes 
and the three numerical second order derivatives $p_1$, $p_2$, and $p_3$, 
and the flexibility of DG methods allows the implantation of this connection 
into the formulation of the proposed mIP-DG methods. 

Like in the finite difference framework, the g-monotonicity (generalized 
monotonicity) and the numerical moment play a central role in the 
proposed mIP-DG framework. The g-monotonicity gives the mIP-DG methods 
the ability to select the mathematically ``correct" solution 
(i.e., the viscosity solution) among all possible solutions, 
and the numerical moment is the catalyst which facilitates 
the g-monotonicity of the proposed mIP-DG methods.
Moreover, the g-monotonicity allows for the possible development of
more efficient (than generic Newton) solvers. The special nonlinearity 
of the algebraic systems can be explored to decouple the equations 
as seen in Algorithm~\ref{solver_alg}. We believe that one of the 
main strengths of the mIP-DG formulation presented in this paper 
lies in the way in which the discretization handles the nonlinearity.  
The discretization takes a nonlinear problem and embeds it into a 
mostly linear system of equations where the nonlinearity has been 
modified to ensure g-monotonicity. The added monotonicity can 
theoretically enlarge the domain of valid initial guesses over 
which a solver will converge. Thus, the weak coupling with
linear equations is only a small penalty for the added
structure in the nonlinearity.

We also remark that the role and benefit of the symmetrizing parameter 
is unclear for nonlinear systems of equations. When $\gamma_0$ is 
sufficiently large, we observe that numerical results seem independent 
of the choice for $\epsilon$.  However, when the penalty constant 
$\gamma_0$ is not sufficiently large, the inclusion of $\epsilon$ 
can be detrimental to the approximation. For small $\gamma_0$, 
the approximation is allowed to have larger jumps occur. When the 
jumps become too large, the effect from having $\epsilon$ present 
becomes exaggerated, and the overall accuracy of the 
approximation begins to suffer beyond just the presence of jumps.  
For elliptic problems, we can see the formation of a boundary layer.  
For dynamic problems, we see the approximation actually diverging 
(almost instantaneously) along the interior of the domain.  
As expected, when $\gamma_0$ increases and becomes sufficiently large, 
these phenomena disappear. Thus, at a numerical level, the presence of 
the symmetrizing constant $\epsilon$ seems important, even though at the 
continuous level for the PDE, the symmetrization terms all become zero. 

Conceptually, the mIP-DG framework presented in this paper can 
be easily extended to the high dimensional fully nonlinear PDE problems, 
the detailed exposition will be given in a forthcoming paper. 
On the other hand, the proposed mIP-DG framework may not work in 
the case when the viscosity solution does not belong 
to $H^1(\Ome)$. In such a case,  a more involved mixed local 
discontinuous Galerkin (mLDG) framework must be invoked. 
We refer the interested reader to \cite{Feng_Lewis12c}
for a detailed exposition.

\bibliographystyle{plain}

\end{document}